\documentclass[12pt]{article}
\usepackage{amsmath,amsthm,amssymb}

\newtheorem{Theorem}{Theorem}
\newtheorem{Lemma}{Lemma}
\newtheorem{Proposition}{Proposition}

\newtheorem{Corollary}{Corollary}

\numberwithin{equation}{section}

\begin{document}
\date{}
\author
{
M.I. Belishev
\thanks
{
St. Petersburg State University,
7--9 Universitetskaya nab., St. Petersburg, 199034 Russia;
m.belishev@spbu.ru.
St. Petersburg Department of V. A. Steklov Institute of Mathematics of the Russian Academy of Sciences,
27 Fontanka, St. Petersburg, 191023 Russia
belishev@pdmi.ras.ru.
The work was supported by the grants RFBR 14-01-00535À and Volkswagen Foundation-2016.
}
,\,
S.A. Simonov
\thanks
{
St. Petersburg Department of V. A. Steklov Institute of Mathematics of the Russian Academy of Sciences,
27 Fontanka, St. Petersburg, 191023 Russia;
St. Petersburg State University,
7--9 Universitetskaya nab., St. Petersburg, 199034 Russia;
St. Petersburg State Technological Institute (Technical University),
26 Moskovsky pr., St. Petersburg, 190013 Russia;
sergey.a.simonov@gmail.com.
The work was supported by the grants RFBR 16-01-00443À and RFBR 16-01-00635À.
}
}

\title{Wave model of the Sturm-Liouville operator on the half-line}

\maketitle

\begin{abstract}
The notion of the {\it wave spectrum} of a semi-bounded symmetric operator was introduced by one of the authors in 2013. The wave spectrum is a topological space determined by the operator in a canonical way. The definition uses a dynamical system associated with the operator: the wave spectrum is constructed from its reachable sets. In the paper we give a description of the wave spectrum of the operator $L_0=-\frac{d^2}{dx^2}+q$ which acts in the space $L_2(0,\infty)$ and has defect indices $(1,1)$. We construct a functional ({\it wave}) model of the operator $L_0^*$ in which the elements of the original $L_2(0,\infty)$ are realized as functions on the wave spectrum. It turns out to be identical to the original $L_0^*$. The latter is fundamental in solving inverse problems: the wave model is determined by their data, which allows for reconstruction of the original.
\end{abstract}

\setcounter{section}{-1}

\section{Introduction}\label{sec Introduction}

The notion of the wave spectrum of a symmetric semi-bounded operator was introduced in \cite{JOT}. The wave spectrum is a topological space determined by the operator in a canonical way. The definition uses a dynamical system associated with the operator: states of the system serve as material for constructing this space. It is constituted of the atoms of the Hilbert lattice of subspaces determined by reachable sets of the system and is endowed with an adequate topology.

The wave spectrum is an invariant of the operator: wave spectra of unitarily-equivalent operators are canonically homeomorphic. At the same time, in important applications the wave spectrum of the operator acting in the space of functions turns out to be homeomorphic to the support of functions which comprise the space. For example, the wave spectrum of the minimal Laplacian on a Riemannian manifold with boundary is essentially identical (isomorphic) to the manifold itself. This fact is used for solving inverse problems. In the problem of reconstruction of a manifold from its boundary data (for instance, the reaction operator) a unitary copy of the Laplacian is extracted from data, and then one can find its wave spectrum. The latter, by construction, is isometric to the manifold and thus gives the solution to the problem (cf. \cite{JOT}, \cite{BD_2}). Such a scheme is a form of the boundary control method (the BC-method). This is an approach to inverse problems which is based on their deep connections with the control theory \cite{BIP97}, \cite{BIP07}.

The notion of the wave spectrum appeared as a result of generalization of the ``experimental material'' gathered from solving particular inverse problems with the BC-method. At some point it became clear that the procedure of solving is equivalent to construction of a certain functional model of the symmetric operator. In this model the elements of the original Hilbert space are realized as functions on the wave spectrum. The outline of this ``wave'' model is given in \cite{JOT}\footnote
{
The idea of this construction can be seen in the early work \cite{KacProblem}. On the heuristic level the wave spectrum was introduced in \cite{BIP07} for solving inverse problems. Formal definition (although with a different name) appeared later, in \cite{Nest}.
}
; its usefulness and effectiveness can be considered as established facts. At the same time, in our opinion, the wave model is also interesting from the theoretical point of view. Its systematic study is our long-term goal.

We consider a particular example: a positive definite Sturm-Liouville operator $L_0=-\frac{d^2}{dx^2}+q$ in $L_2(0,\infty)$ which has defect indices $(1,1)$. We construct the wave model of the operator $L_0^*$. As we proceed, we describe the elements of the general construction and in parallel clarify how exactly they look in our case. At some point, realizing the elements of the original $L_2(0,\infty)$ as complex-valued functions on the wave spectrum, we use specifics of the Sturm-Liouville operator. In the general case the realization is more complex: corresponding functions map to linear spaces of rather abstract nature \cite{JOT}. The above-mentioned specifics allows for complete investigation of the model. The wave spectrum turns out to be identical to the half-line $[0,\infty)$, and the model operator is related with the original $L_0^*$ by a simple gauge transform. As a consequence, the potential $q$ is easily recovered, which determines the original operator.

We dedicate this work to the memory of Vladimir Savel'evich Buslaev, a wonderful person, an excellent mathematician, one of our Teachers.

\section{Dynamical system with boundary control}\label{sec DSBC}

\subsection{The operator $L_0$}\label{ssec Operator $L_0$}

Let us describe the class of operators for which the definitions that we give make sense. Let $\cal H$ be a (separable) Hilbert space and $L_0$ be an operator in $\cal H$. We suppose that
\begin{enumerate}
\item
    $L_0$ is closed and densely defined: $\overline{{\rm Dom\,}L_0} = \cal H$,
\item
    $L_0$ is positive definite: for some $\varkappa >0$ and for every $y \in {\rm Dom\,}L_0$ one has $(L_0 y,y)\geqslant \varkappa \|y\|^2$,
\item
    $L_0$ has non-zero defect indices $n_\pm=\dim {\rm Ker\,} L_0^* \leqslant \infty$.
\end{enumerate}
From the third assumption it follows that $L_0$ is unbounded. Let us denote by $L$ the Friedrichs extension of $L_0$: $L_0 \subset L \subset L_0^*,\,\,L^*= L$ and $(L y,y)\geqslant\varkappa \|y\|^2$ for every $y \in {\rm Dom\,}L$ (cf. \cite{BirSol}). The inverse operator $L^{-1}$ is bounded and is defined on the whole space $\cal H$.
\bigskip

\noindent$\bullet$\,\,\,Everywhere below $H^s$ are Sobolev classes; ${{\mathbb R}_+}:=(0, \infty),\, {\bar{\mathbb R}_+} :=[0,\infty)$. ``Smooth'' always means ``$C^\infty$-smooth''.

In the case of Sturm-Liouville operator we have ${\cal H} = L_2({\mathbb R}_+)$. The operator is $L_0: {\cal H} \to {\cal H}$,
\begin{align}
\notag & {\rm Dom\,} L_0 = \left\{y \in {\cal H}\cap H^2_{\rm loc}(\bar{\mathbb R}_+)\,|\,\,y(0)=y'(0)=0;\,\,-y''+qy \in {\cal H}\right\},
\\
\label{SL first}& L_0y:=-y''+qy\,,
\end{align}
where $q=q(x)$ is a real-valued function {(\it the potential)} such that
\begin{equation}\label{q conditions}
\begin{array}{l}
(1)\,\,\,\,\,\, q\in C^{\infty}(\bar{\mathbb R}_+),
\\
(2) \,\,\,\,\text{ the limit point case takes place},
\\
(3)\,\,\,\, \text{ the operator } L_0 \text{ is positive definite.}
\end{array}
\end{equation}
In such a case the problem
\begin{equation}\label{phi}
-\phi''+q\phi=0\,,\,\,\,x>0\,; \qquad \phi(0)=1\,, \qquad \phi \in
L_2({{\mathbb R}_+})
\end{equation}
has a unique solution $\phi(x)$ which is a smooth function. It is known that, firstly,
\begin{align}
\notag & {\rm Dom\,} L_0^* = \left\{y \in {\cal H}\cap H^2_{\rm
loc}(\bar{\mathbb
R}_+)\,|\,-y''+qy \in {\cal H}\right\},\\
\notag & L_0^*y:=-y''+qy\,,\\
\label{L_0^*} & {\rm Ker\,}L_0^*=\{c \phi\,|\,\,c\in\mathbb C\},
\end{align}
and the defect indices of $L_0$ are $n_{\pm}={\rm dim\,}{\rm Ker\,}L_0^*=1$\footnote
{
This holds, for instance, if $q(x)> - cx^2$ with some $c>0$: cf. \cite{Naimark}, Chapter VII.26, Theorem 6.
}
; secondly, the Friedrichs extension of the operator $L_0$ is $L: {\cal H} \to {\cal H}$,
\begin{align}
\notag & {\rm Dom\,} L = \left\{y \in {\cal H}\cap H^2_{\rm
loc}(\bar{\mathbb
R}_+)\,|\,\,y(0)=0;\,\,-y''+qy \in {\cal H}\right\},\\
\label{L}& Ly:=-y''+qy\,.
\end{align}
We should add that smoothness of $q$ simplifies considerations, however, all the main results can be extended to the case $q \in C_{\rm loc}\left({\bar{\mathbb R}_+}\right)$ at least.

\subsection{Green's system}\label{ssec Green system}

The following definitions are close to the ones used in the classical work of A. N. Kochubei \cite{Koch} (see also \cite{Ryzh}, \cite{MMM}).
\smallskip

Let $\cal H$ and $\cal B$ be Hilbert spaces, $A: {\cal H} \to\cal H$ and $\Gamma_i: {\cal H} \to {\cal B}\,\,\,(i=1,2)$ be two operators for which the following conditions hold:
$$
    \overline{{\rm Dom\,}A}={\cal H},\quad{\rm Dom\,}
    \Gamma_i \supset {\rm Dom\,}A,\quad\overline{{\rm
    Ran\,}\Gamma_1 +{\rm
    Ran\,}\Gamma_2} = \cal B.
$$
The set ${\mathfrak G}=\{{\cal H}, {\cal B}; A, \Gamma_1, \Gamma_2\}$ is called the {\it Green's system}, if its elements are related by the {\it Green's formula}
\begin{equation}\label{Green Formula}
(Au,v)_{\cal H}-(u,Av)_{\cal H}=(\Gamma_1u, \Gamma_2 v)_{\cal B} -
(\Gamma_2u, \Gamma_1 v)_{\cal B}
\end{equation}
for every $u, v \in {\rm Dom\,} A$. The space $\cal H$ is called {\it inner}, $\cal B$ is {\it the space of boundary values}, $A$ is {\it the basic operator}, $\Gamma_{1,2}$ are {\it boundary operators}.

\subsection{The system ${\mathfrak G}_{L_0}$}\label{ssec System Gr_L_0}

One can relate to the operator $L_0$ satisfying the conditions 1.-- 3. of Section
\ref{ssec Operator $L_0$} a Green's system in a canonical way. Let
$$
    {\cal K}\,:=\,{\rm Ker\,} L_0^*\,.
$$
Denote by $P_K$ the orthogonal projection in $\cal H$ on ${\cal K}$, by $\mathbb O$ and $\mathbb I$ denote zero and identity operators. Define also the following operators:
\begin{equation}\label{Ga1Ga2}
    \Gamma_1\,:=\,L^{-1}L_0^*-{\mathbb I}\,, \quad \Gamma_2\,:=\,P_K L_0^*.
\end{equation}
As it was shown in \cite{BD_2} (Lemma 1), the set ${\mathfrak G}_{L_0}:=\{{\cal H}, {\cal K};\,L_0^*, \Gamma_1, \Gamma_2\}$ forms a Green's system. In the same work action of the boundary operators $\Gamma_i$ was described in terms of the {\it Vishik's decomposition}, which has the following form:
\begin{equation}\label{Vishik general}
    {\rm Dom\,}L_0^*={\rm Dom\,}L_0 \overset{.}+L^{-1}{\cal K}\overset{.}+{\cal K}
\end{equation}
(the sums are direct). If one applies this decomposition to an arbitrary $y\in {\rm Dom} L_0^*$,
\begin{equation}\label{Vishik concrete}
    y=y_0+L^{-1}g_y + h_y, \quad y_0 \in {\rm Dom}\, L_0, \quad g_y, h_y \in {\mathcal K},
\end{equation}
then the boundary operators (\ref{Ga1Ga2}) act by the rule
\begin{equation}\label{**}
    \Gamma_1 y = -h_y, \quad \Gamma_2 y = g_y
\end{equation}
(cf. \cite{BD_2}, Section 2.3).
\smallskip

\noindent$\bullet$\,\,\,For the Sturm-Liouville operator (\ref{SL first}) we have ${\cal K}=\{c \phi\,|\,\,c\in\mathbb C\}$. Let $\eta:=L^{-1}\phi$ and note that $\eta(0)=0$ owing to (\ref{L}), while $\eta'(0)\not=0$. Indeed, if one assumes $\eta(0)=\eta'(0)=0$, one then has $\eta \in {\rm Dom\,} L_0$ and
$L_0\eta=L_0L^{-1}\phi=\phi$ which, owing to
$$
0=(L_0^*\phi,\eta)=(\phi,L_0\eta)=(\phi,\phi)\not=0,
$$
leads to a contradiction.

It is easy to check that in our case the representations (\ref{Vishik concrete}) and (\ref{**}) take the form
\begin{align}
\notag &
y=\left\{y-y(0)\phi-\left[\frac{y'(0)-y(0)\phi'(0)}{\eta'(0)}\right]\eta\right\}+
\left[\frac{y'(0)-y(0)\phi'(0)}{\eta'(0)}\right]\eta+y(0)\phi\,;\\
\label{Ga1Ga2 SL} & \Gamma_1y = -\,y(0)\phi\,, \qquad \Gamma_2
y=\left[\frac{y'(0)-y(0)\phi'(0)}{\eta'(0)}\right]\phi
\end{align}
(recall that $\phi(0)=1$).

In this way to the operator (\ref{SL first}) canonically corresponds a Green's system with the inner space $L_2({\mathbb R}_+)$, the basic operator (\ref{L_0^*}), the boundary space $\{c \phi\,|\,\,c\in\mathbb C\}$ and the boundary operators (\ref{Ga1Ga2 SL}).

\subsection{The system $\alpha_{L_0}$}\label{ssec System alpha}

The system ${\mathfrak G}_{L_0}$ that corresponds to the operator $L_0$ determines in turn a dynamical system
\begin{align}
\label{alpha1} & u_{tt}+L_0^*u = 0,  && t>0,\\
\label{alpha2} & u|_{t=0}=u_t|_{t=0}=0, && \\
\label{alpha3}& \Gamma_1 u = h\,, && t\geqslant 0,
\end{align}
where $h=h(t)$ is the {\it boundary control} (${\cal K}$-valued function of time), $u=u^h(t)$ is the solution (${\cal H}$-valued function of time). In the control theory, $u^h(\cdot)$ is called the trajectory, $u^h(t)$ is the state of the system at the moment $t$. Having applications in mind we call $u^h$ the {\it wave}. The system (\ref{alpha1})--(\ref{alpha3}) is determined by the operator $L_0$, and we denote it by ${\alpha_{L_0}}$.

Recall that $L$ is the Friedrichs extension of the operator $L_0$. By $L^{\frac{1}{2}}$ let us denote the positive square root of $L$. Assume that the control $h$ is smooth and vanishes near $t=0$. Denote by
\begin{equation}\label{linal cal M}
{\cal M}:=\{h \in C^\infty\left([0,\infty); {\cal K}\right)\,|\,\,{\rm supp\,}h \subset (0,\infty)\}
\end{equation}
a linear span of such controls. As is shown in \cite{BD_2}, for $h \in \cal M$ the problem (\ref{alpha1})--(\ref{alpha3}) has a unique classical solution $u^h \in C^\infty\left([0,\infty); {\cal H}\right)$. It vanishes near $t=0$ and admits the following representation:
\begin{equation}
    \label{weak solution u^f}
    u^h(t)\,=\,-h(t)+\int_0^t L^{-\frac{1}{2}}\,\sin\left[(t-s)L^{\frac{1}{2}}\right]\,h_{tt}(s)\,ds\,,
    \qquad t \geqslant 0\,.
\end{equation}
For controls from the class $\{h\,|\,\,h, h_{tt} \in L_2^{\rm loc}\left([0,\infty); {\cal K}\right),\,\,h(0)=h_t(0)=0\}$ the (generalized) solution is defined as the right-hand side of (\ref{weak solution u^f}). To distinguish generalized solutions from classical let us call the latter {\it smooth waves}. In what follows they will play the role of a certain structure in $\cal H$.

\noindent$\bullet$\,\,\,In the case of the Sturm-Liouville operator we have ${\cal K}=\{c \phi\,|\,\,c\in\mathbb C\}$, and the condition (\ref{alpha3}) takes the form $\Gamma_1 u =h(t)=f(t)\phi,\,\,t\geqslant 0\,$, with some complex-valued function $f$. Hence the system (\ref{alpha1})--(\ref{alpha3}) is equivalent to the following initial boundary value problem:
\begin{align}
\label{alpha1*} & u_{tt}-u_{xx}+qu = 0,  && x>0,\,\,t>0,\\
\label{alpha2*} & u|_{t=0}=u_t|_{t=0}=0, && x \geqslant 0,\\
\label{alpha3*} & u|_{x=0} = f\,, && t\geqslant 0\,.
\end{align}
Let us define an analog of (\ref{linal cal M}), the lineal
\begin{equation}\label{linal cal M SL}
{\dot{\cal M}}:=\{f \in C^\infty\left[0,\infty\right)\,|\,\,{\rm
supp\,}f \subset (0,\infty)\}
\end{equation}
of smooth controls that vanish near $t=0$. For $f \in \cal M$ the problem (\ref{alpha1*})--(\ref{alpha3*}) has a unique classical solution $u=u^f(x,t)$, which is smooth in both variables. For this solution the following representation holds:
\begin{equation}\label{u^f repres}
u^f(x,t)=f(t-x)+\int_x^t w(x,s)\,f(t-s)\,ds\,, \qquad x \geqslant
0, \,\,t \geqslant 0\,,
\end{equation}
where both summands on the right-hand side are considered to vanish for $x>t$. The function $w$ is defined for $0\leqslant x\leqslant t$ and is smooth; it is simply related to the classical Riemann function of the equation (\ref{alpha1*}). This representation is used to define the (generalized) solution corresponding to controls
$f \in L_2^{\rm loc}[0,\infty)$: such a solution is {\it defined} as the right-hand side of (\ref{u^f  repres}).

Solution $u^f(\cdot,t)$ considered as an $L_2({\mathbb R}_+)$-valued function of time is the trajectory of the system $\alpha$ that corresponds to the operator (\ref{SL first}). This case is specific in that the map $f \mapsto u^f$ is continuous: it is easy to see that $u^f \in C_{\rm loc}\left([0,\infty); L_2({\mathbb R}_+)\right)$.

\subsection{Controllability}\label{ssec Controllability}

Let us return to the system $\alpha_{L_0}$ in the general case.  Fix $t=T\geqslant 0$; the set of states
\begin{equation}\label{U^t}
{\cal U}^T_{L_0}\,:=\,\{u^h(T)\,|\,\,h \in \cal M\}
\end{equation}
is called {\it reachable} (at the time $T$). It is easy to see that ${\cal U}_{L_0}^T$ grows with $T$. Define also
\begin{equation}\label{U}
    {\cal U}_{L_0}\,:=\,
  \bigcup_{T > 0}{\cal U}^T_{L_0}\,,\quad {\cal D}_{L_0}\,:=\,{\cal H} \ominus \overline{\cal
  U}_{L_0},
\end{equation}
the {\it total} reachable set and the {\it defect subspace}. The lineal of smooth waves is invariant under $L_0^*$. Indeed, for $u=u^h(T)\in {\cal U}_{L_0}^T$ one has
\begin{align*}
&L_0^*u^h(T)
\overset{(\ref{alpha1})}=-u^h_{tt}(T)\overset{(\ref{weak solution u^f})}=-u^{h_{tt}}(T)\in
{\cal U}_{L_0}^T,
\\
&u^h(T) \overset{(\ref{alpha2})}=
J^2[u^h_{tt}](T)\overset{(\ref{weak solution u^f})}=u^{J^2h}_{tt}(T)
\overset{(\ref{alpha1})}=-L_0^*u^{J^2h}(T)\in L_0^*{\cal U}_{L_0}^T\,,
\end{align*}
where $J:=\int_0^t$ is integration in time. From this we conclude that $L_0^*{\cal U}_{L_0}= {\cal U}_{L_0}$.

Let us remark in advance that the functional model of the operator $L_0^*$ which we are constructing is in fact the model of its {\it wave part} $L_0^*|_{{\cal U}_{L_0}}$. Related to this is the following question left unanswered: let $\overline{\cal U}_{L_0}=\cal H$, i.e., the part $L_0^*|_{{\cal U}_{L_0}}$ is densely defined; is it then true that its closure coincides with $L_0^*$? In all the examples we know the answer is positive.

A system ${\alpha_{L_0}}$ is called {\it controllable}, if
\begin{equation}\label{contr}
    \overline{{\cal U}}_{L_0}\,=\,{\cal H} \qquad \left({\cal D}_{L_0}=\{0\}\right)\,.
\end{equation}
Let us formulate a criterion of controllability, which was established in \cite{BD_2}.

Recall the definitions. An operator $A$ is said to induce a self-adjoint operator in a (non-trivial) subspace  ${\cal L}\subset \cal H$, if $\overline{{\cal L} \cap {\rm Dom\,}A}=\cal L$, $A \left[{\cal L} \cap {\rm Dom\,}A\right]\,\subset \cal L$ and the operator $A|_{{\cal L} \cap {\rm Dom\,}A}$ is self-adjoint in $\cal L$. The operator $A$ is called {\it completely non-self-adjoint}, if there does not exist a subspace in $\cal H$ where $A$ induces a self-adjoint operator. As is shown in \cite{BD_2} (Theorem 1), the system ${\alpha_{L_0}}$ is controllable, if and only if  $L_0$ is a non-self-adjoint operator.
\medskip

\noindent$\bullet$\,\,\,In the case of the Sturm-Liouville operator \eqref{SL first} let us show that the system ${\alpha}$ is controllable. Let $C^\infty_{\rm fin}\left({\bar{\mathbb R}_+}\right)$ be the set of smooth functions with bounded support. Denote
 $$
C^\infty_T\left({\bar{\mathbb R}_+}\right)\,:=\,\{y \in C^\infty\left({\bar{\mathbb R}_+}\right)\,|\,\,{\rm supp\,}y\subset [0,T)\}\,.
 $$
Obviously, $C^\infty_{\rm fin}\left({\bar{\mathbb R}_+}\right)=\cup_{T>0}C^\infty_T\left({\bar{\mathbb R}_+}\right)$.

\begin{Lemma}\label{Lemma controll SL}
Let the operator $L_0$ have the form \eqref{SL first} and the potential $q$ satisfy the conditions \eqref{q conditions}. Then the following relations hold:
\begin{equation}\label{UT SL}
{\cal U}^T\,=\,C^\infty_T\left({\bar{\mathbb R}_+}\right),\quad
T>0\,; \qquad {\cal U}\,=\, C^\infty_{\rm fin}\left({\bar{\mathbb
R}_+}\right),
\end{equation}
where ${\cal U}^T$ and ${\cal U}$ are defined in \eqref{U^t} and \eqref{U}.
\end{Lemma}

$\square$
\,\,\, Fix $T>0$ and pick an $f \in \dot{\cal M}$. According to (\ref{u^f repres}) we have
$$
u^f(x,T)\,=\,f(T-x)+\int_x^T w(x,s)\,f(T-s)\,ds\,, \qquad x\geqslant 0\,.
 $$
From this one can easily see that the wave $u^f(\cdot,T)$ is a smooth function which vanishes near $x=T$. Hence, the left-hand side of the first of the equalities (\ref{UT SL}) is a subset of the right-hand side.

Pick a $y$ from the right-hand side. Let us find $f=f(t)|_{0\leqslant t\leqslant T}$ from the Volterra integral equation of the second kind
 $$
f(T-x)+\int_x^T w(x,s)\,f(T-s)\,ds\,=\,y(x), \qquad 0 \leqslant x
\leqslant T\,.
 $$
It is easy to see that $f$ is a smooth function which vanishes near $t=0$. Continue $f$ to the interval $(T,\infty)$ in an arbitrary way that preserves its smoothness. By construction we have: $f \in \cal M$ è $y=u^f(\cdot,T) \in {\cal U}^T$. Therefore, the right-hand side of the first equality of (\ref{UT SL}) is a subset of the left-hand side. Thus the equality is established.

The second equality follows from the first.
\quad
$\blacksquare$
\bigskip

As a consequence, we have controllability: $\overline{{\cal U}}=\overline{C^\infty_{\rm fin}\left({\bar{\mathbb R}_+}\right)}=L_2({\mathbb R}_+)$. From controllability it follows that the operator (\ref{SL first}) is non-self-adjoint. This fact can be also proved directly, without using dynamics. Besides that, using (\ref{UT SL}) it is not difficult to show that the closure of the wave part $L_0^*|_{\cal U}$ of the operator (\ref{L_0^*}) coincides with $L_0^*$ itself.

\section{The wave spectrum}\label{sec Wave spectrum}

\subsection{The wave isotony}\label{ssec Isotonicity}

Let  $\cal P$ and $\cal Q$  be partially ordered sets. A map $i:\,{\cal P}\to \cal Q$ is called {\it isotone} (order-preserving), if $p \preceq p'$  implies $i(p) \preceq i(p')$ \cite{Birkhoff}.

We call a family of isotone maps $\{i^t\}_{t \geqslant 0}$ such that $p \preceq p'$ and $t\leqslant t'$ imply $i^t(p) \preceq i^{t'}(p')$ an {\it isotony}. In another formulation, an isotony is an isotone map of the set ${\cal P}\times [0,\infty)$ (with the natural order on it) into $\cal Q$.
\smallskip

{\it Lattice} is a partially ordered set in which every two elements $x,y$ have the least upper bound $x \vee y$ and the greatest lower bound $x \wedge y$ (cf. \cite{Birkhoff}). We will deal with concrete lattices endowed with additional structures: complements, topology, etc.

Let ${{\mathfrak L}(\cal H)}$ be the lattice of (closed) subspaces of $\cal H$ with the partial order $\subseteq$. Its is easy to check that ${\cal A}\vee {\cal B}=\overline{\{a+b\,|\,\,a \in{\cal A},\, b \in \cal B\}}$ è $\cal A\wedge\cal B={\cal A} \cap{\cal B}$. The lattice ${{\mathfrak L}(\cal H)}$ is also a lattice with the least, $\{0\}$ , and the greatest, $\cal H$, elements, and with complements, ${\cal A}^\bot={\cal H}\ominus \cal A$ (since ${\cal A}^{\bot}\vee{\cal A}=\cal H$, ${\cal A}^{\bot}\wedge{\cal A}=\{0\}$). By $P_{\cal A}$ we denote the (orthogonal) projection in $\cal H$ on $\cal A$. The topology on ${{\mathfrak L}(\cal H)}$ is defined by the strong operator convergence of projections: ${\cal A}_j \to {\cal A}$, if $P_{{\cal A}_j}\overset{s}\to P_{\cal A}$ as $j\rightarrow\infty$\footnote
{
The strong operator topology is not first countable, and thus it cannot be described in terms of convergence of sequences. However, its restriction to the set of orthogonal projections (as well as to any subset of $\mathcal B(\mathcal H)$ bounded in the operator norm) is first countable and even metrizable
\cite{Kim-2006}.
}
.

\medskip

By \emph{lattice isotony} ${{\mathfrak L}(\cal H)}$ we call an isotony $I=\{I^t\}_{t \geqslant 0}: {{\mathfrak L}({\cal H})\times [0,\infty)} \to {{\mathfrak L}(\cal H)}$ with the following additional properties: $I^0\,=\,{\rm id}$ and $I^t\{0\}\,=\,\{0\}$. In what follows we deal only with such $I$.
\smallskip

To every operator $L_0$ from the class defined in Section \ref{ssec Operator $L_0$} a \emph{wave isotony} $I_{L_0}$ corresponds in the following way. Consider the dynamical system
\begin{align}
\label{dual 1} & v_{tt}+L v = g\,,  && t>0,\\
\label{dual 2} & v|_{t=0}=v_t|_{t=0}=0\,,
\end{align}
where $g$ is an $\cal H$-valued function of time. If $g\in C^\infty\left([0,\infty); {\cal H}\right)$ vanishes near $t=0$, then the problem has a unique classical solution $v=v^g(t)$, for which the Duhamel's representation holds: \begin{align}\label{v^g general}
    v^g(t)=\int_0^t
    L^{-\frac{1}{2}}\,\sin\left[(t-s)L^{\frac{1}{2}}\right]g(s)\,ds\,,
    \qquad t \geqslant 0
\end{align}
(cf. \cite{BirSol}). For $g \in L_2^{\rm loc}\left([0,\infty); {\cal H}\right)$ the (generalized) solution is defined as the right-hand side (\ref{v^g general}).

Fix a subspace $\cal G \in {{\mathfrak L}(\cal H)}$ and consider $\cal G$-valued controls. Corresponding reachable sets of the system (\ref{dual 1})--(\ref{dual 2}) are
\begin{equation}\label{V^t}
    {\cal V}^t_{\cal G}:=\left\{v^g(t)\,|\,\,g \in L_2^{\rm loc}\left([0,\infty); {\cal G}\right)\right\}.
\end{equation}
It is clear that ${\cal V}^t_{\cal G}$ grows with ${\cal G}$ and $t$. Define the family of maps ${I_{L_0}}=\{I_{L_0}^t\}_{t \geqslant 0}$:
\begin{equation}\label{I^t}
    I_{L_0}^0:={\rm id}\,; \qquad I_{L_0}^t{\cal G}:=\overline{{\cal V}^t_{\cal G}}\,,\,\,\,\,\,t>0\, .
\end{equation}

\begin{Proposition}\label{Prop 1}
The family ${I_{L_0}}$ is an isotony of the lattice ${\mathfrak L}(\cal H)$.
\end{Proposition}

The proof can be found in \cite{JOT}. Note that $I_{L_0}$ is determined not by $L_0$, but by its Friedrichs extension $L$. It is clear that the wave isotony can be correctly defined for every self-adjoint operator which is semi-bounded from below. In applications, the problem (\ref{dual 1})--(\ref{dual 2}) describes propagation of waves excited by the sources $g$, so that the initial subspace $\cal G$ is extended by the waves $v^g$.
\medskip

\noindent$\bullet$\,\,\,Let us discuss properties of the wave isotony for the Sturm-Liouville operator  (\ref{SL first}). For a subset $E$ of the half-line denote by $E^r$ its $r$-neighborhood:
$$
E^r:=\{x \in \bar{\mathbb R}_+\,|\,\,{\rm dist\,}(x,E):=\underset{e \in E}{\rm inf\,}|x-e|<r\}\,, \qquad r>0\,.
$$
Let $\Delta_{a,b}$ be one of the intervals $(a,b),\,[a,b),\,(a,b],\,[a,b] \,\,\,(0\leqslant a<b\leqslant\infty)$. Then, obviously,
\begin{equation}\label{metric nbh Delta}
\Delta_{a,b}^r\,=\,\begin{cases}
          (a-r,\,b+r)\,, & a\geqslant r,\\
          [0,\,b+r)\,, & a<r.
                   \end{cases}
\end{equation}
If $x_0 \in\bar{\mathbb R}_+$ is a point, then
\begin{equation}\label{metric nbh x0}
\{x_0\}^r\,=\,\begin{cases}
          (x_0-r,\,x_0+r)\,, & x_0\geqslant r,\\
          [0,\,x_0+r)\,, & x_0<r.
                   \end{cases}
\end{equation}
For a measurable set $E \subset {\mathbb R}_+$ let us denote$L_2(E):=\{y \in L_2({\mathbb R}_+)\,|\,\,y|_{CE}=0\}$, where $CE:={\mathbb R}_+ \setminus E$.

\begin{Lemma}\label{Lemma Basic}
Under  conditions of Lemma \ref{Lemma controll SL}, for every $0\leqslant a<b\leqslant\infty$ and $T>0$ the following relation holds:
\begin{equation}\label{L2(Delta) to L-2(Delta^r)}
I^TL_2(\Delta_{a,b})=L_2(\Delta^T_{a,b})\,,
\end{equation}
where $I^T$ is defined by \eqref{I^t}.
\end{Lemma}

$\square$
\,\,\,{\bf 1.}\,\,\,In our case the system (\ref{dual 1})--(\ref{dual 2}) is equivalent to the initial boundary problem
\begin{align}
\label{dual 1 SL} & v_{tt}-v_{xx}+qv = g\,,  && x \in {\mathbb R}_+,\,\,t>0,\\
\label{dual 2 SL} & v|_{t=0}=v_t|_{t=0}=0\,, && x \in
{\bar{\mathbb R}_+},\\
\label{dual 3 SL} & v|_{x=0}=0\,, && t \geqslant 0
\end{align}
with the right-hand side $g = g(x,t)$ such that $g(\cdot,t)\in {\cal G}$ for every $t\geqslant 0$. The condition (\ref{dual 3 SL}) follows from $v(\cdot,t)\in {\rm Dom\,} L$, according to (\ref{L}). The problem (\ref{dual 1 SL})--(\ref{dual 3 SL}) is well-posed for every $g \in L^{\rm loc}_2\left({\mathbb R}_+\times [0,\infty)\right)$; for $g\in C^\infty_0\left({\mathbb R}_+\times (0,\infty)\right)$ its solution $v=v^g(x,t)$ is a classical one. Besides that, owing to finiteness of the domain of influence for the (hyperbolic) equation (\ref{dual 1 SL}), if ${\rm supp\,} g(\cdot,t) \subset\overline{\Delta_{a,b}}$  for every $t>0$, then ${\rm supp\,}v^g(\cdot,T) \subset \overline{\Delta^T_{a,b}}$ for every $T>0$.

Fix $T>0$ and choose $0 \leqslant a<b\leqslant\infty$. Consider ${\cal G}=L_2(\Delta_{a,b})$ and take $g \in L_2^{\rm loc}\left(\Delta_{a,b} \times [0,\infty)\right)$. From (\ref{v^g general}) it follows that the map $L_2[0,T]\ni g|_{[0,T]}\mapsto v^g(\cdot,T)\in L_2({\mathbb R}_+)$ is continuous. Besides that, ${\rm supp\,} v^g(\cdot,T) \subset\overline{\Delta^T_{a,b}}$ holds. Therefore the following inclusion takes place: ${\cal
V}^T_{\cal G} \subset L_2(\Delta^T_{a,b})$. Let us show that ${\cal V}^T_{\cal G}$ is dense in $L_2(\Delta^T_{a,b})$.

{\bf 2.}\,\,\, Consider an auxiliary problem
\begin{align}
\label{d-dual 1 SL} & w_{tt}-w_{xx}+qw = 0\,,  && x \in {\mathbb R}_+,\,\,0<t<T,\\
\label{d-dual 2 SL} & w|_{t=T}=0,\,\,w_t|_{t=T}=y\,, && x \in
{\bar{\mathbb R}_+},\\
\label{d-dual 3 SL} & w|_{x=0}=0\,, && 0\leqslant t \leqslant T\,.
\end{align}
It is well-posed for every $y \in L_2^{\rm loc}\left({\mathbb R}_+\right)$, and for $y \in C^\infty_0({\mathbb R}_+)$ its solution $w=w^y(x,t)$ is classical. By finiteness of the domain of influence, if ${\rm supp\,} y \subset \overline{\Delta_{a,b}}$, then ${\rm supp\,} w^y(\cdot,t) \subset \overline{\Delta^{T-t}_{a,b}}$ for every $t>0$.

Let us establish relations between solutions of the problems (\ref{dual 1 SL})--(\ref{dual 3 SL}) and (\ref{d-dual 1 SL})--(\ref{d-dual 3 SL}).  Let $g \in C^\infty_0\left({\mathbb R}_+\times(0,\infty)\right)$ and $y \in C^\infty_0({\mathbb R}_+)$, so that the solutions of both problems are classical (smooth). By finiteness of the domain of influence, for every $t>0$ the functions $v^g(\cdot,t)$ and $w^y(\cdot,t)$ have compact supports in ${\bar{\mathbb R}_+}$. This fact justifies the following calculation.

Integrating by parts we have:
\begin{align*}
& \int_{{\mathbb R}_+\times[0,T]}g w^y\,dx dt\overset{(\ref{dual 1 SL})}=\int_{{\mathbb R}_+\times[0,T]}[v^g_{tt}-v^g_{xx}+qv^g]w^y\,dx dt=
\\
&=\int_{{\mathbb R}_+}[v^g_t w^y-v^g w^y_t]\big|^{t=T}_{t=0}dx-\,\int^T_0[v^g_x w^y-v^gw^y_x]\big|^{x=\infty}_{x=0}dt-
\\
& -\int_{{\mathbb R}_+\times[0,T]}v^g[w^y_{tt}-w^y_{xx}+qw^y]\,dxdt \overset{(\ref{dual 2 SL})-(\ref{d-dual 3 SL})}=-\int_{{\mathbb R}_+}v^g(\cdot,T)\,y\, dx\,,
\end{align*}
and
\begin{equation}\label{duality relation}
\int_{{\mathbb R}_+\times[0,T]}g w^y\,dx dt=-\int_{{\mathbb R}_+}v^g(\cdot,T)\,y\, dx\,.
\end{equation}
Since $C^\infty_0({\mathbb R}_+)$ is dense in  $L_2({\mathbb R}_+)$, the last equality holds for every $y\in L_2({\mathbb R}_+)$.

Let $g \in C^\infty_0\left({\mathbb R}_+\times(0,\infty)\right)$ and $y \in L_2({\mathbb R}_+)$; suppose additionally that ${\rm supp\,}g(\cdot,t)\subset \Delta_{a,b}$ for every $t>0$. It follows that ${\rm
supp\,}v^g(\cdot,T)\subset\overline{\Delta^T_{a,b}}$, and the equality (\ref{duality relation}) takes the form
\begin{equation}\label{duality relation final}
\int_{\Delta_{a,b}\times[0,T]}g w^y\,dx dt=-\left(v^g(\cdot,T),y\right)_{L^2(\Delta^T_{a,b})}\,.
\end{equation}

{\bf 3.}\,\,\,Return to the question of density of ${\cal V}^T_{\cal G}$ in $L_2(\Delta^T_{a,b})$. Let $y\in
L_2(\Delta^T_{a,b})\ominus \overline{{\cal V}^T_{\cal G}}$ and show that $y=0$.

From the choice of $y$, the right-hand side of (\ref{duality relation final}) is $0$. Since $g$ is arbitrary, we conclude that
\begin{equation}\label{***}
w^y\,=\,0 \quad \text{in}\,\,\,\Delta_{a,b}\times [0,T]\,.
\end{equation}
Let us continue the solution $w^y$ to the times $T\leqslant t\leqslant 2T$ by oddness:
$$
w^y(\cdot,t)\,:=\,
\begin{cases}
w^y(\cdot,t), & 0\leqslant t< T,
\\
w^y(\cdot, 2T-t), & T\leqslant t\leqslant 2T.
\end{cases}
$$
The continuation solves the following problem:
\begin{align}
\label{dd-dual 1 SL} & w_{tt}-w_{xx}+qw = 0\,,  && x \in {\mathbb R}_+,\,\,0<t<2T,\\
\label{dd-dual 2 SL} & w|_{t=T}=0,\,\,w_t|_{t=T}=y\,, && x \in
\bar{\mathbb R}_+,\\
\label{dd-dual 3 SL} & w|_{x=0}\,=\,0\,, && 0\leqslant t \leqslant
2T\,.
\end{align}
One only needs to check that the continued function $w^y$ satisfies the equation (\ref{dd-dual 1 SL}). This is easy, because the odd continuation does not result in jumps of $w^y$ and $w^y_t$ at $t=T$. Owing to (\ref{***}) we have:
\begin{equation}\label{****}
w^y\,=\,0 \quad \text{in}\,\,\,\Delta_{a,b}\times [0,2T]\,.
\end{equation}

The solution of the equation (\ref{dd-dual 1 SL}) with the property (\ref{****}) can be continued by zero from $\Delta_{a,b}\times [0,2T]$ to the wider domain
$$
\Omega^{2T}_{a,b}:= \left\{(x,t)\,|\,\,\max\{0,|t-T|+(a-T)\}<x<b+T-|t-T|\right\}\,,
$$
which is bounded by the corresponding characteristic lines of the equation (\ref{dd-dual 1 SL}). Indeed, fix the point $(x_0,t_0)$ for which $b<x_0<b+T-|t_0-T|$. This point belongs to $\Omega^{2T}_{a,b}$, but not to $\Delta_{a,b}\times [0,2T]$. Take a small $\delta>0$ such that the characteristic triangle (the influence con-e) $\{(x,t)\,|\,\,b-\delta<x<x_0-|t-t_0|\}$ is contained in $\Omega^{2T}_{a,b}$. By finiteness of the domain of influence, the value $w^y(x_0,t_0)$ is determined by the Cauchy data $w^y,\,w^y_x$ on the (vertical) base of the cone, which is a subset of the line $x=b-\delta$. This base is contained inside $\Delta_{a,b}\times [0,2T]$, and owing to (\ref{****}) this data is zero. Hence $w^y(x_0,t_0)=0$. One can consider the points $(x_0,t_0)\in\Omega^{2T}_{a,b}$ with $x_0<a$ (if any) in an analogous fashion.

Summing up, $w^y=w^y_t=0$ everywhere in $\Omega^{2T}_{a,b}$. In particular, for $t=T$ according to (\ref{dd-dual 2 SL}) we have $y=w^y_t(x,T)=0$ for every $x \in \Delta^T_{a,b}$. Therefore from $y \in L_2(\Delta^T_{a,b})\ominus \overline{{\cal V}^T_{\cal G}}$ it follows that $y=0$. Thus the relation $\overline{{\cal V}^T_{\cal G}}=L_2(\Delta^T_{a,b})$ is established, i.e., (\ref{L2(Delta) to L-2(Delta^r)}) is proved. \qquad
 $\blacksquare$

\subsection{Lattices, atoms, the wave spectrum}\label{ssec Lattices and atoms}

{\it Lattice} in ${{\mathfrak L}(\cal H)}$ is a subset invariant under all the operations in ${{\mathfrak L}(\cal H)}$ that were defined in the beginning of Section \ref{ssec Isotonicity}. Every lattice necessarily contains $\{0\}$ and $\cal H$.

Let ${\mathfrak M}\subset {\mathfrak L}(\cal H)$ be a family of subspaces. By ${\mathfrak L}_{\mathfrak M}$ we denote the minimal ${\mathfrak L}(\cal H)$ that contains ${\mathfrak M}$. It consists of all the subspaces of the form $\vee_{1\leqslant k \leqslant n}\cap_{1\leqslant l \leqslant m}{\cal A}_{kl}$ where for each subspace ${\cal A}_{kl}$ either ${\cal A}_{kl}\in \mathfrak M$ or ${\cal A}^\bot_{kl}\in \mathfrak M$ holds (cf. \cite{BirSol}).
\smallskip

Let $I$ be an isotony of the lattice ${\mathfrak L}(\cal H)$. The family ${\mathfrak M}\subset{\mathfrak L}(\cal H)$ is {\it invariant} under $I$, if $I{\mathfrak M}:=\{I^t {\cal M}\,|\,\,t\geqslant 0,\,\,{\cal M}\in{\mathfrak M}\}={\mathfrak M}$. For every ${\mathfrak M}\subset {\mathfrak L}(\cal H)$ there exists a minimal lattice ${\mathfrak L}^I_{\mathfrak M}$ which contains $\mathfrak M$ and is invariant under $I$. It is easy to check that it has the following constructive description. Define the operation $\sigma$ on subsets of the lattice ${\mathfrak L}(\cal H)$ by the rule $\sigma({\mathfrak M}):=I{\mathfrak L}_{\mathfrak M}$. Then ${\mathfrak L}^I_{\mathfrak M}=\cup_{j\geqslant 1}\sigma^j({\mathfrak M})$.

Let ${\cal F}\left([0,\infty); {\mathfrak L}(\cal H)\right)$ be the set of ${{\mathfrak L}(\cal H)}$-valued functions of $t$. It is a lattice with point-wise defined partial order, operations and convergence:
\begin{align*}
    & \{f \leqslant g \} \Longleftrightarrow \{f(t) \subseteq g(t),\,\,t \geqslant 0\}, \quad (f \vee g)(t):=f(t)\vee g(t),
    \\
    & (f \wedge g)(t):=f(t)\cap g(t),\,\,(f^\bot)(t):=(f(t))^\bot,\,\,(\lim f_j)(t):=\lim (f_j(t))\,.
\end{align*}
The least and the greatest elements of this lattice are the functions $0_{\cal F}$ and $1_{\cal F}$ identically equal to $\{0\}$ and $\cal H$, respectively. If $\mathfrak L\subset{\mathfrak L}(\cal H)$ is a lattice, then the set ${\cal F}\left([0,\infty); {\mathfrak L}\right)$, which consists of $\mathfrak L$-valued functions, is also a lattice. If $\mathfrak L$ is invariant under the isotony (of the lattice) $I$, then the set of motonically growing functions
$$
{\cal F}_I\left([0,\infty); {\mathfrak L}\right):=\{f(t)=I^t{\cal L}\,|\,\,{\cal L}\in \mathfrak L\}
$$
is contained in ${\cal F}\left([0,\infty); {\mathfrak L}\right)$. In what follows its completion $\overline{{\cal F}_I\left([0,\infty); {\mathfrak L}\right)}\subset {\cal F}\left([0,\infty); {\mathfrak L}(\cal H)\right)$, a set obtained by adding to ${\cal F}_I\left([0,\infty); {\mathfrak L}\right)$ the limits of all converging sequences in ${\cal F}_I\left([0,\infty); {\mathfrak L}\right)$, plays an important role.
\smallskip

In the most general setting, let $\cal P$ be a partially ordered set with the least element $0$. The element $\omega \in \cal P$ is called an {\it atom}, if
$\omega \not= 0$ and from $0\not=\omega'\preceq \omega$ follows $\omega'=\omega$ (cf. \cite{Birkhoff}). By ${\rm At}\, \cal P$ we denote the set of all atoms of $\cal P$.
\smallskip

Consider the system $\alpha_{L_0}$. Recall that its reachable sets are defined by (\ref{U^t}); let ${\mathfrak U}_{L_0}:=\{\,\overline{{\cal U}^T_{L_0}}\,\}_{T \geqslant0}\subset {\mathfrak L}(\cal H)$ be the family of reachable subspaces (closures of reachable sets). The family ${\mathfrak U}_{L_0}$ and the wave isotony $I_{L_0}$ are determined by the operator $L_0$. As a consequence, the operator determines the (minimal) lattice ${\mathfrak L}_{L_0}:={\mathfrak L}^{I_{L_0}}_{{\mathfrak U}_{L_0}}$ which contains all reachable subspaces and is invariant under $I_{L_0}$. The lattice and the isotony determine the set of functions ${\cal F}_{I_{L_0}}\!\left([0,\infty); {\mathfrak L}_{L_0}\right)$. Thus there is canonical correspondence between the operator $L_0$ and the set of atoms
$$
\Omega_{L_0}\,:=\,{\rm At}\,\overline{{\cal F}_{I_{L_0}}\left([0,\infty); {\mathfrak L}_{L_0}\right)}\,.
$$
This set is called the {\it wave spectrum} of the operator $L_0$ and is the main object of interest in this paper.

Certain additional assumptions about the operator $L_0$ provide that ${\Omega_{L_0}}\not=\emptyset$, cf. \cite{JOT}. There exist operators with the wave spectrum which consists of a single point. In the general case the question of whether $\Omega_{L_0}$ is non-empty remains open.
\medskip

\noindent$\bullet$\,\,\,Let us look for the case of the Sturm-Liouville operator $L_0$ given by (\ref{SL first}) at the objects defined above. According to (\ref{UT SL}), reachable subspaces of the corresponding system $\alpha$ are ${\mathfrak U}=\{L_2(\Delta_{0,T})\}_{T\geqslant 0}$. Action of the wave isotony $I$ on subspaces $L_2(\Delta_{a,b})$ is described by Lemma \ref{Lemma Basic}.

Let us call the set $E \subset {\bar{\mathbb R}_+}$ {\it elementary}, if $E=\cup_{j=1}^{n(E)}\Delta_{a_j,b_j}$, where $0\leqslant a_1<b_1<a_2<b_2<...<a_{n(E)}<b_{n(E)}\leqslant \infty$. Let ${\cal E}$ be the family of all elementary sets. Obviously, the metric extension $E\mapsto E^T=\{x \in {\bar{\mathbb R}_+}\,|\,{\rm dist\,}(x,E)<T\}$ (cf. (\ref{metric nbh Delta}) and  (\ref{metric nbh x0})) maps elementary sets to elementary sets. We will call subspaces $L_2(E)$ with $E \in \cal E$ {\it elementary}. The family of such subspaces forms the lattice ${\mathfrak L}_{\mathcal E}\subset {\mathfrak L}({\cal H})$.

\begin{Lemma}\label{Lemma Basic 2}
Under conditions of Lemma \ref{Lemma controll SL} one has $I^TL_2(E)=L_2(E^T)$ for every $E\in\cal E$.
\end{Lemma}

$\square$
\,\,\, The set $E$ can be written as $E=\cup_{j=1}^{n(E)}\Delta_{a_j,b_j}$. Owing to isotonicity of $I_{L_0}$, we obtain from Lemma \ref{Lemma Basic} that $L_2(\Delta^T_{a_j,b_j})=I^TL_2(\Delta_{a_j,b_j})\subseteq I^TL_2(E)$ for every $j=1,2,...,n(E)$, and therefore $L_2(E^T)\subseteq I^TL_2(E)$. The same argument as in the first part of the proof of Lemma \ref{Lemma Basic} leads to the inclusion ${\cal V}^T_{L_2(E)}\subseteq L_2(E^T)$ and, hence, $I^T(L_2(E))\subseteq L_2(E^T)$ holds.
\quad
$\blacksquare$

It is important that $\cal E$ contains only unions of intervals of positive length (non-degenerate): for a degenerate interval (as $E=\{x\}$, $x\in\mathbb R_+$) the equality $I^TL_2(E)=L_2(E^T)$ obviously fails.

It follows from the construction of the minimal lattice $\mathfrak L_{\mathfrak M}$ for $\mathfrak M\subset\mathfrak L(\mathcal H)$ given above that $\mathfrak L_{\mathfrak U}=\mathfrak L_{\cal E}$. From Lemma \ref{Lemma Basic 2} it follows directly that the lattice $\mathfrak L_{\cal E}$ is invariant under the wave isotony $I$, and therefore the equality $\mathfrak L=\mathfrak L_{\mathfrak U}^{I}=\mathfrak L_{\cal E}$ holds.

Below, $m$ is the Lebesgue's measure on $\bar{\mathbb R}_+$, $A\triangle B:=[A\setminus B]\cup[B\setminus A]$ is the symmetric difference of the sets $A$ and $B$. Let ${\rm Leb}(\mathbb R_+)$ be the $\sigma$-algebra of Lebesgue-measurable subsets of the half-line $\bar{\mathbb R}_+$.

\begin{Lemma}\label{Lemma convergence in L(H)}
Let $\{E_n\}_{n=1}^{\infty}$ be a sequence of sets from ${\rm Leb}(\mathbb R_+)$, $E\in{\rm Leb}(\mathbb R_+)$. Then $L_2(E_n)\overset{\mathfrak L(\mathcal H)}\longrightarrow L_2(E)$ as $n\rightarrow\infty$, if and only if for every $L>0$ one has $m((E_n\triangle E)\cap(0,L))\rightarrow0$ as $n\rightarrow\infty$.
\end{Lemma}

$\square$
\,\,\, If $L_2(E_n)\overset{\mathfrak L(\mathcal H)}\longrightarrow L_2(E)$ as $n\rightarrow\infty$, then the projections on these subspaces converge in the strong sense: $P_{L_2(E_n)}\overset{s}\rightarrow P_{L_2(E)}$. Then for every $L>0$ it is true that $P_{L_2(E_n)}\chi_{(0,L)}\overset{L_2(\mathbb R_+)}\longrightarrow P_{L_2(E)}\chi_{(0,L)}$, where $\chi_{(0,L)}$ is the characteristic function of the interval $(0,L)$. This means that $\int_0^L|\chi_{E_n}(x)-\chi_E(x)|^2dx=m((E_n\triangle E)\cap(0,L))\rightarrow0$ as $n\rightarrow\infty$.

Let us now suppose that for every $L>0$ one has $m((E_n\triangle E)\cap(0,L))\rightarrow0$ as $n\rightarrow\infty$. Then $P_{L_2(E_n)}\chi_{(0,L)}\overset{L_2(\mathbb R_+)}\longrightarrow P_{L_2(E)}\chi_{(0,L)}$. Linear span of the set $\{\chi_{(0,L)}|L>0\}$ is dense in $L_2(\mathbb R_+)$. By the Banach-Steinhaus theorem we have $P_{L_2(E_n)}\overset{s}\rightarrow P_{L_2(E)}$ as $n\rightarrow\infty$, and this by definition means that
$L_2(E_n)\overset{\mathfrak L(\mathcal H)}\longrightarrow L_2(E)$ as $n\rightarrow\infty$. \quad $\blacksquare$
\bigskip

In accordance with the lemma let us define a convergence in ${\rm Leb}(\mathbb R_+)$ as follows: a sequence $\{E_n\}_{n=1}^{\infty}\subset{\rm Leb}(\mathbb R_+)$ \textit{converges} to $E\in{\rm Leb}(\mathbb R_+)$, if for every $L>0$ one has $m((E_n\triangle E)\cap(0,L))\rightarrow0$ as $n\rightarrow\infty$. Denote $\mathfrak L_{\rm Leb(\mathbb R_+)}:=\{L_2(E)|E\in{\rm Leb}(\mathbb R_+)\}$.

\begin{Lemma}\label{Lemma closure of L-L-0 is L-Leb}
Under conditions of Lemma \ref{Lemma controll SL} the following equality holds: $\overline{\mathfrak L}=\mathfrak L_{\rm Leb(\mathbb R_+)}$, where the closure is taken in the sense of convergence in $\mathfrak L(\mathcal H)$.
\end{Lemma}

$\square$
\,\,\,Let us first show that $\overline{\mathfrak L}\subseteq\mathfrak L_{\rm Leb(\mathbb R_+)}$. Let $\{L_2(E_n)\}_{n=1}^{\infty}$ be a sequence from $\mathfrak L$ which is convergent in the sense of the topology on $\mathfrak L(\mathcal H)$. We need to show that its limit $A\in\mathfrak L(\mathcal H)$ belongs to $\mathfrak L_{\rm Leb(\mathbb R_+)}$. From existence of the limit it follows that for every $L>0$ the sequence of functions $P_{L_2(E_n)}\chi_{(0,L)}=\chi_{E_n\cap(0,L)}$ converges in $L_2(\mathbb R_+)$ and, therefore, is fundamental. The following equalities hold:
\begin{multline}\label{1}
\|\chi_{E_n\cap(0,L)}-\chi_{E_m\cap(0,L)}\|_{L_2(\mathbb R_+)}=\int_0^L|\chi_{E_n}(x)-\chi_{E_m}(x)|^2dx
\\
=m((E_n\cap(0,L))\triangle(E_m\cap(0,L))).
\end{multline}
For every $L>0$ the function $\rho_L(F,G)=m(F\triangle G)$, $F,G\in{\rm Leb}(0,L)$, is a pseudometric on ${\rm Leb}(0,L)$. The equivalence relation $F\sim F'$, if $\rho_L(F,F')=0$, determines the set of equivalence classes ${\rm Leb}^{\sim}(0,L)$, which is a complete metric space \cite{KolFom}. From \eqref{1} it follows that the sequence $\{E_n\cap(0,L)\}_{n=1}^{\infty}$ is fundamental in ${\rm Leb}(0,L)$. Thus the sequence of equivalence classes
$\{(E_n\cap(0,L))^{\sim}\}_{n=1}^{\infty}$ converges to some equivalence class $E^{\sim}(L)$. For $L_2>L_1$ the intersection of every representative of the equivalence class $E^{\sim}(L_2)$ with the interval $(0,L_1)$ belongs to the equivalence class $E^{\sim}(L_1)$. Hence there exists a set $E\in{\rm Leb}(\mathbb R_+)$ such that for every $L>0$ one has $E\cap(0,L)\in E^{\sim}(L)$. This exactly means that $E_n\rightarrow E$ as $n\rightarrow\infty$ in the sense of the above definition, which by Lemma \ref{Lemma convergence in L(H)} means convergence $L_2(E_n)\overset{\mathfrak L(H)}\longrightarrow L_2(E)$ as $n\rightarrow\infty$. Therefore $A=L_2(E)\in\mathfrak L_{\rm Leb(\mathbb R_+)}$.

Let us now show that $\mathfrak L$ is dense in $\mathfrak L_{\rm Leb(\mathbb R_+)}$. This is equivalent to density of $\cal E$ in ${\rm Leb}(\mathbb R_+)$ in the sense of the above definition. It is enough to show that for every $L>0$ the set $\{E\in{\cal E}|E\subseteq(0,L)\}$ is dense in ${\rm Leb}(0,L)$ with respect to the pseudometric $\rho_L$. Every measurable subset of $(0,L)$ can be approximated in $\rho_L$ by open subsets of $(0,L)$, therefore open subsets are dense in ${\rm Leb}(0,L)$. Every open bounded set on the real line is an at most countable union of non-intersecting open intervals such that the sequence of their lengths is summable. Thus every open subset of $(0,L)$ can be approximated in pseudometric $\rho_L$ by a finite union of intervals, and therefore in such a way one can approximate every measurable subset of $(0,L)$. This means that $\cal E$ is dense in ${\rm Leb}(\mathbb R_+)$. \quad
$\blacksquare$
\bigskip

For $x\geqslant 0$ denote
\begin{equation}\label{omega x}
\omega_x(t):=L_2(\{x\}^t).
\end{equation}
$\omega_x(t)$ is a monotone function of $t$ with values in $\mathfrak L(\mathcal H)$, i.e., is an element of ${\cal F}\left([0,\infty); {\mathfrak L}\right)\subset{\cal F}\left([0,\infty); {\mathfrak L}(\cal H)\right)$.

\begin{Lemma}\label{Lemma omega-x are atoms}
Under conditions of Lemma \ref{Lemma controll SL} for every $x\in\bar{\mathbb R}_+$ one has $\omega_x\in\Omega$.
\end{Lemma}

$\square$
\,\,\,
Let us first show that $\omega_x\in\overline{{\cal F}_{I}\left([0,\infty); {\mathfrak L}\right)}$. If $x=0$, then for every $t\geqslant0$
\begin{equation*}
\{0\}^t=[0,t)=\lim\limits_{n\rightarrow\infty}\left(0,t+\frac1n\right)
=\lim\limits_{n\rightarrow\infty}\left(0,\frac1n\right)^t,
\end{equation*}
where both limits are in the sense of convergence in ${\rm Leb}(\mathbb R_+)$. If $x>0$, then, owing to \eqref{metric nbh x0}, for every $t\geqslant0$
\begin{multline*}
\{x\}^t=\lim\limits_{n\rightarrow\infty}\left(\max\left\{0,x\left(1-\frac1{n}\right)-t\right\},
x\left(1+\frac1{n}\right)+t\right)
\\
=\lim\limits_{n\rightarrow\infty}\left(x\left(1-\frac1{n}\right),x\left(1+\frac1{n}\right)\right)^t.
\end{multline*}
Hence, by Lemma \ref{Lemma convergence in L(H)}, the corresponding subspaces converge as $n\rightarrow\infty$ in $\mathfrak L(\mathcal H)$ for every $t\geqslant 0$, which means that $\omega_x\in\overline{{\cal F}_{I}\left([0,\infty);{\mathfrak L}\right)}$.

Let us now show that $\omega_x$ is an atom of the lattice $\overline{{\cal F}_{I}\left([0,\infty); {\mathfrak L}\right)}$. Suppose there exists a non-zero element $\omega\in\overline{{\cal F}_{I}\left([0,\infty); {\mathfrak L}\right)}$ such that $\omega\leqslant \omega_x$. For every $t\geqslant 0$ one can write $\omega(t)=L_2(E(t))$ with some measurable set $E(t)\subseteq\{x\}^t$. Since $\omega\in\overline{{\cal F}_{I}\left([0,\infty); {\mathfrak L}\right)}$, there exists a sequence $\{E_n\}_{n=1}^{\infty}$ in $\cal E$ such that for every $t\geqslant 0$ one has $I^t(L_2(E_n))\overset{\mathfrak L(\mathcal H)}\longrightarrow L_2(E(t))$ or, by Lemma \ref{Lemma convergence in L(H)}, $E_n^t\overset{{\rm Leb}(\mathbb R_+)}\longrightarrow E(t)$ as $n\rightarrow\infty$. For $\delta\in(0,1)$ denote $F_n(\delta):=E_n\cap(x(1-\delta),x(1+\delta))$. For every $t\geqslant 0$ we have the inclusion $E(t)\subseteq\{x\}^t$. Thus for every $L>0$ one has $m(((F_n(\delta))^t\triangle E(t))\cap(0,L))\leqslant m((E_n^t\triangle E(t))\cap(0,L))\rightarrow0$, which means $(F_n(\delta))^t\overset{{\rm Leb}(\mathbb R_+)}\longrightarrow E(t)$ as $n\rightarrow\infty$.

For every set $E\in\cal E$ the derivative $\frac{d(m(E^t))}{dt}$ is equal to the number of positive edges of non-intersecting open intervals comprising the set $E^t$. For the set $F_n(\delta)$ this number is at least two for  $t<x(1-\delta)$ and at least one for
$t\geqslant x(1-\delta)$, so $m((F_n(\delta))^t)\geqslant m(F_n(\delta))+\min\{2t,x(1-\delta)+t\}$. Taking the limit, we get $m(E(t))\geqslant \min\{2t,x(1-\delta)+t\}$ for arbitrarily small $\delta>0$. This means that $m(E(t))\geqslant\min\{2t,x+t\}=m(\{x\}^t)$ and, since $E(t)\subseteq\{x\}^t$, one has $m(E(t)\triangle\{x\}^t)=0$ for every $t\geqslant 0$. Therefore every non-zero element $\omega\in\overline{{\cal F}_{I}\left([0,\infty); {\mathfrak L}\right)}$ such that $\omega\leqslant\omega_x$ should coincide with $\omega_x$. Thus $\omega_x$ is an atom of the lattice $\overline{{\cal F}_{I}\left([0,\infty); {\mathfrak L}\right)}$. \quad
$\blacksquare$
\bigskip

The following result characterizes the wave spectrum of the operator (\ref{SL first}).

\begin{Theorem}\label{Th 1}
Let $L_0$ be the operator given by \eqref{SL first} with $q$ satisfying the conditions \eqref{q conditions}. Then the set $\Omega$ is bijective to the half-line ${\bar{\mathbb R}_+}:$
$$
\Omega=\{\omega_x\,|\,\,x\in{\bar{\mathbb R}_+}\},
$$
where the elements $\omega_x$ are defined by \eqref{omega x}.
\end{Theorem}

$\square$
\,\,\, By Lemma \ref{Lemma omega-x are atoms}, $\{\omega_x|x\geqslant 0\}\subseteq\Omega$. To prove the inverse inclusion take an atom $\omega\in\Omega$. For every $t\geqslant 0$ the subspace $\omega(t)$ has the form $\omega(t)=L_2(E(t))$ where $E(t)$ is some measurable set. Hence there exists a sequence $\{E_n\}_{n=1}^{\infty}$ in $\cal E$ such that $E_n^t\overset{{\rm Leb}(\mathbb R_+)}\longrightarrow E(t)$ as $n\rightarrow\infty$.

For every $L>0$ one has $E_n^t\cap(0,L+t)\overset{{\rm Leb}(\mathbb R_+)}\longrightarrow E(t)\cap(0,L+t)$ as $n\rightarrow\infty$. There exist $t_0\geqslant 0$ and $L_0>0$ such that $m(E(t_0)\cap(0,L_0+t_0))>0$. Hence there exists $N_0$ such that for every $n\geqslant N_0$ one has $m(E_n^{t_0}\cap(0,L_0+t_0))>0$. Since
\begin{equation*}
E_n^{t_0}\cap(0,L_0+t_0)=(E_n\cap(0,L_0+2t_0))^{t_0}\cap(0,L_0+t_0),
\end{equation*}
for every $n\geqslant N_0$ one has $m((E_n\cap(0,L_0+2t_0))^{t_0})>0$. Denote $L_1:=L_0+2t_0$ and $F_n:=E_n\cap(0,L_1)$. Since $E_n\in\cal E$, $F_n$ does not contain degenerate intervals; since $m(F_n^{t_0})>0$, $F_n\neq\emptyset$. Therefore  $F_n\in\cal E$ and $m(F_n)>0$.

For every $L>0$ and $n,m\in\mathbb N$ we have:
\begin{multline*}
m((F_n^t\triangle F_m^t)\cap(0,L))=m((((E_n\cap(0,L_1))^t)\triangle(E_m\cap(0,L_1))^t)\cap(0,L))
\\
\leqslant m((E_n^t\triangle E_m^t)\cap(0,L)).
\end{multline*}
The sequence $\{E_n^t\cap(0,L)\}_{n=1}^{\infty}$ is fundamental in pseudometric $\rho_L$ for every $L>0$, so the sequence $\{F_n^t\}_{n=1}^{\infty}$ is also fundamental in $\rho_L$. Thus it has a limit which we denote by $F(t)$ (the set $F(t)$ is defined not uniquely, but up to a set of measure zero). Since for every $n\geqslant N_0$ we have $F(t)\backslash E(t)\subseteq(F(t)\backslash F_n^t)\cup(E_n^t\backslash E(t))$, and so
$$
m(F(t)\backslash E(t))\leqslant m((F(t)\backslash F_n^t))+m((E_n^t\backslash E(t))),
$$
we obtain $m(F(t)\backslash E(t))=0$. $\omega$ is an atom, therefore $\omega(t)=L_2(F(t))$.

For every $n\geqslant N_0$ one has: $F_n\in\cal E$, $F_n\subseteq(0,L_1)$, $m(F_n)>0$. For $t>L_1$ from Lemma \ref{Lemma Basic 2} we have $F_n^t=(0,\sup F_n+t)$. Existence of the limit $F_n^t$ as $n\rightarrow\infty$ in ${\rm Leb}(\mathbb R_+)$ means that for every $\varepsilon>0$ there exists $N_1(\varepsilon)$ such that for every $n,m\geqslant N_1(\varepsilon)$ one has
\begin{multline*}
    m(F_n^t\triangle F_m^t)
    =m([\min\{\sup F_n,\sup F_m\}+t,\max\{\sup F_n,\sup F_m\}+t))
    \\
    =|\sup F_n-\sup F_m|<\varepsilon.
\end{multline*}
Therefore $\{\sup F_n\}_{n=1}^{\infty}$ is a fundamental sequence of positive numbers. Denote its limit by $L_2$. For every $\varepsilon>0$ there exists $N_2(\varepsilon)$ such that for every $n\geqslant N_2(\varepsilon)$ one has $\sup F_n\in(L_2-\varepsilon,L_2+\varepsilon)$. Then for $t>\varepsilon$ the following inclusion takes place:
$$
(\max\{0,L_2+\varepsilon-t\},L_2-\varepsilon+t)\subseteq
F_n^t.
$$
Since $F_n^t\overset{{\rm Leb}(\mathbb R_+)}\longrightarrow F(t)$, one has
\begin{multline*}
m((\max\{0,L_2+\varepsilon-t\},L_2-\varepsilon+t)\backslash F(t))
\\
\leqslant
m((\max\{0,L_2+\varepsilon-t\},L_2-\varepsilon+t)\backslash F_n^t)+m(F_n^t\backslash F(t))
\\
=m(F_n^t\backslash F(t))\rightarrow0\text{ as }n\rightarrow\infty
\end{multline*}
and, consequently, $m((\max\{0,L_2+\varepsilon-t\},L_2-\varepsilon+t)\backslash F(t))=0$. This holds for every $\varepsilon\in(0,t)$, so
$$
m((\max\{0,L_2-t\},L_2+t)\backslash F(t))=m(\{L_2\}^t\backslash
F(t))=0
$$
for every $t\geqslant 0$. $\omega$ is an atom, therefore $\omega=\omega_{L_2}$.
\quad$
\blacksquare$

\subsection{The space $\Omega_{L_0}$}\label{ssec Space Omega L0}

Let us return to the case of a general $L_0$. The wave spectrum, if it is not empty, can be naturally endowed with certain structures.
\medskip

\noindent{\bf Topology.}
\,\,\, By definition, atoms are ${\mathfrak L}{({\cal H})}$-valued functions of time. Fix an atom $\omega\in{\Omega_{L_0}}$: $\omega=\omega(t)$, $t\geqslant 0$. The set
\begin{equation}\label{balls}
    B_r[\omega]:=\{\omega' \in {\Omega_{L_0}}\,|\,\exists\,t>0:\,\{0\}\not=
    \omega'(t)\subseteq \omega(r)\} \qquad (r>0)
\end{equation}
is called a {\it ball}, $\omega$ and $r$ are its center an radius.

\begin{Proposition}\label{Prop 2}
The system of balls $\{B_r[\omega]\,|\,\,\omega \in {\Omega_{L_0}},\,\,r>0\}$ is a base of some topology on $\Omega_{L_0}$.
\end{Proposition}

The proof is given in \cite{BD_2}; it checks characteristic properties of a base. Thus the wave spectrum becomes a topological space.

There exist other natural topologies on $\Omega_{L_0}$. Relations between them are yet to be revealed, cf. \cite{JOT}. The ball topology now seems to us the most relevant. However, its general properties (Hausdorffness, metrizability, etc.) are not studied.
\smallskip

\noindent{\bf Metric.}
\,\,\, Under additional assumptions about atoms one can introduce a metric on $\Omega_{L_0}$. To each atom $\omega\in{\Omega_{L_0}}$: $\omega=\omega(t)$, $t\geqslant 0$ corresponds a positive operator in $\cal H$
\begin{equation}\label{eikonal}
\tau_\omega\,:=\,\int_{[0,\infty)} t\,dP_{\omega(t)}\,,
\end{equation}
where $P_{\omega(t)}$ is the projection on the subspace $\omega(t)\subseteq \cal H$. We call it an {\it eikonal\,}; this term is motivated by applications, cf. \cite{JOT, BD_2}. Introduce the distance
\begin{equation}\label{eikonal dist}
\tau: \Omega_{L_0}\times \Omega_{L_0}\to [0,\infty)_,\quad
\tau(\omega,\omega'):=\|\tau_\omega-\tau_{\omega'}\|\,.
\end{equation}
Below we will see that this definition can be correct even in the case of unbounded $\tau_\omega$. However, in the general case one cannot exclude a pathologic situation of $\tau=\infty$. How the ball topology is related to the topology that corresponds to the metric (\ref{eikonal dist}), is also an open question.
\smallskip

\noindent{\bf The boundary.}
\,\,\,Let us return to the system ${\alpha_{L_0}}$ and the family $\{\overline{{\cal U}^T_{L_0}}\}_{T \geqslant 0\,}\subset{\mathfrak L}(\cal H)$ of its reachable subspaces. The set of atoms
\begin{equation}\label{boundary}
    \partial {\Omega_{L_0}}\,:=\,\left\{\omega \in {\Omega_{L_0}}\,|\,\,\,\omega(t)
    \subseteq \overline{{\cal U}^t_{L_0}},\,\forall\,t>0\right\}
\end{equation}
is called the {\it boundary} of the wave spectrum. Whether $\partial
{\Omega_{L_0}}$ is always non-empty, is an open question.
\medskip

\noindent$\bullet$
\,\,\,Let $L_0$ be the Sturm-Liouville operator (\ref{SL first}). By
\begin{equation}\label{bijection beta}
\beta: \,{\bar{\mathbb R}_+}\ni x\, \mapsto \,\omega_x \in\Omega
\end{equation}
we denote the canonical bijection established by Theorem \ref{Th 1}. Below, ${\rm dist}(x,x')=|x-x'|$ is the standard distance in ${\bar{\mathbb R}_+}$.

\begin{Lemma}\label{Lemma eikonal}
Let $\omega \in \Omega$ be an atom. The eikonal corresponding to it is the unbounded self-adjoint operator  $\tau_\omega$ with the domain
$$
{\rm Dom\,} \tau_\omega=\left\{y \in L_2({\mathbb R}_+)\,|\,\,\int_0^\infty(1+x)^2|y(x)|^2\,dx<\infty\right\}.
$$
Its action is multiplication by the distance:
\begin{equation}\label{+}
\left(\tau_\omega y\right)(x)\,=\,{\rm dist}(x,x_\omega)\,y(x)\,,
\quad x \in {\bar{\mathbb R}_+}\,,
\end{equation}
where  $x_\omega=\beta^{-1}(\omega)$.
\end{Lemma}

$\square$
\,\,\,Pick a function $y \in L_2({\mathbb R}_+)$ with a compact support.

By Theorem \ref{Th 1} we have: $\omega(t)=L_2(\{x_\omega\}^t),\,\,t\geqslant 0$. Hence the operator $P_{\omega(t)}$, which projects $L_2({\mathbb R}_+)$ on $L_2(\{x_\omega\}^t)$, acts by cutting functions to the neigh\-bor\-hood $\{x_\omega\}^t$:
\begin{equation}\label{++}
\left(P_{\omega(t)}y\right)(x)\,=\,
\begin{cases}
y(x)\,, & {\rm dist\,}(x, x_\omega)<t\\
0\,,   & {\rm dist\,}(x, x_\omega)>t
\end{cases}\,.
\end{equation}

Let $T>0$ be such that ${\rm supp\,}y\subset \{x_\omega\}^T$. Take a partition $\Xi:=\{t_i\}_{i=0}^N:\,\,t_0<t_1<...<t_N$ of the interval $\overline{\{x_{\omega}\}^T}$ and the points $\tilde t_i \in[t_{i-1},t_i]$. The value $r_\Xi:=\max_{1\leqslant i \leqslant
N}(t_i-t_{i-1})$ is the rank of the partition. From the definition of integral we have in (\ref{eikonal}):
 $$
\tau_\omega y\,=\,\lim \limits_{r_\Xi\to 0}\sum
\limits_{i=1}^N\tilde t_i\,\Delta_i P_{\omega(t)} y\,,
 $$
where $\Delta_i P_{\omega(t)}:=P_{\omega(t_i)}-P_{\omega(t_{i-1})}$ and convergence is in $L_2({\mathbb R}_+)$ norm. Since $\Delta_iP_{\omega(t)}\Delta_jP_{\omega(t)}=\mathbb O$ for $i\not=j$, the summands are pairwise orthogonal. By (\ref{++}), they equal
\begin{equation*}
\left(\tilde t_i\,\Delta_i P_{\omega(t)}y\right)(x)\,=\,
         \begin{cases}
     \tilde t_iy(x)\,, & {\rm dist\,}(x, x_\omega)<t_i-t_{i-1}\\
     0\,,   & {\rm dist\,}(x, x_\omega)>t_i-t_{i-1}
         \end{cases}\,.
\end{equation*}
In the first line we have $\tilde t_i={\rm dist\,}(x,x_\omega)+O(r_\Xi)$ uniformly in $x \in {\rm supp\,}y$ and $i=1,...\,N$. Using the equality $y=\sum\limits_{i=1}^N \Delta_iP_{\omega(t)} y$ and orthogonality of summands, we get
\begin{align*}
& \left\|{\rm dist}(\cdot, x_\omega)\,y - \sum
\limits_{i=1}^N\tilde t_i\,\Delta_i P_{\omega(t)}
y\right\|^2=\left\|\sum \limits_{i=1}^N\left[{\rm dist}(\cdot,
x_\omega)-\tilde
t_i\right]\,\Delta_i P_{\omega(t)} y\right\|^2=\\
&=\sum \limits_{i=1}^N\left[{\rm dist}(\cdot, x_\omega)-\tilde
t_i\right]^2\left\|\Delta_i P_{\omega(t)} y\right\|^2=O(r_\Xi^2)\sum
\limits_{i=1}^N\left\|\Delta_i P_{\omega(t)}
y\right\|^2=O(r_\Xi^2)\|y\|^2\,.
\end{align*}
Passing to the limit as $r_\Xi \to 0$, we arrive at (\ref{+}).

Closure extends $\tau_\omega$ from functions with finite support to the natural domain ${\rm Dom\,} \tau_\omega=\left\{y \in L_2({\mathbb R}_+)\,|\,\,\int_0^\infty[1+{\rm dist}(x,x_\omega)]^2|y(x)|^2\,dx<\infty\right\}$. Con\-di\-tions $y\in{\rm Dom\,}\tau_\omega$ and $\int_0^\infty(1+x)^2|y(x)|^2\,dx<\infty$ are obviously equivalent.\quad
$\blacksquare$

\begin{Corollary}\label{Cor 1}
The function (\ref{eikonal dist}) defines a metric on $\Omega$; moreover,
\begin{equation}\label{isometry}
\tau(\omega,\omega')\,=\,{\rm dist}(x_{\omega},x_{\omega'})\,,
\qquad \omega,\omega' \in \Omega\,.
\end{equation}
\end{Corollary}

Indeed, we have
\begin{align*}
\tau(\omega,\omega')=\|\tau_{\omega}-\tau_{\omega'}\|\overset{(\ref{+})}=\underset{x
\in {\bar{\mathbb R}_+}}{\rm sup}|{\rm dist}(x, x_{\omega})-{\rm
dist}(x, x_{\omega'})|={\rm dist}(x_{\omega},x_{\omega'})\,.
\end{align*}
From (\ref{isometry}) we conclude that the bijection $\beta$ is an isometry from ${\bar{\mathbb R}_+}$ (with the metric $\rm dist$) to $\Omega$ (with the metric $\tau$). The following facts can be seen from this.

\begin{Proposition}\label{Prop 3}
The balls (\ref{balls}) are identical with the balls corresponding to the $\tau$-metric:
$B_r[\omega]=\{\omega'\in\Omega\,|\,\,\tau(\omega,\omega')<r\}$, so that the ball topology on $\Omega$ coincides with the metric topology. There exists the unique measure $\nu$ on $\Omega$ such that
\begin{equation}\label{mes nu}
\nu\left(B_r[\omega]\right)=m(\{x_\omega\}^r)=r+{\rm min}\{r,x_\omega\}\,.
\end{equation}
The boundary $\partial \Omega$ of the wave spectrum consists of the single atom $\omega_0=\beta(0)$. The function $\tau:\Omega\to[0,\infty),\,\tau(\omega):=\tau(\omega,\omega_0)=x_\omega$ is a global coordinate on $\Omega$.
\end{Proposition}

We omit here simple check of these facts. Let us only note that $\partial\Omega=\{\omega_0\}$ follows from the definition of the boundary (\ref{boundary}) and the equalities
$$
\overline{{\cal U}^T}\overset{(\ref{UT SL})}=\overline{C^\infty_T\left({\bar{\mathbb R}_+}\right)}=L_2(\Delta_{0,T})=L_2(\{0\}^T)=\omega_0(T) \qquad (T \geqslant0)\,,
$$
the last of which is established by Theorem \ref{Th 1}.

We write ${\bar{\mathbb R}_+}[\cdot]$ to specify the variable that we consider. The coor\-di\-na\-ti\-za\-tion
\begin{equation}\label{Coord spectrum}
\Omega\ni \omega\, \mapsto\, \tau(\omega)\in[0,\infty)=:{\bar{\mathbb R}_+}[\tau]
\end{equation}
makes the wave spectrum an isometric copy of the original half-line ${\bar{\mathbb R}_+}[x]$. Summing up these considerations, we see that the wave spectrum of the Sturm-Liouville operator on the half-line with the defect indices $(1,1)$ is in fact identical to the half-line itself.

\section{The wave model}\label{sec Model}

\subsection{The spaces $\tilde{\cal H}$ and ${\cal H}^{\rm w}$}\label{ssec Space tilde cal H}

Let $L_0$ be an operator in $\cal H$  with a non-zero wave spectrum. The wave model is devised to realize elements $y \in \cal H$ as functions $\tilde y(\cdot)$ on $\Omega_{L_0}$ with values in ``natural'' auxiliary spaces. A universal way to map $y \mapsto \tilde y(\cdot)$ was proposed in \cite{JOT} and is described below.
\smallskip

\noindent{\bf Germs.}
\,\, Fix $\omega\in\Omega_{L_0}:\,\,\omega=\{\omega(t)\}_{t \geqslant 0}$. Recall that $P_{\omega(t)}$ is the projection on $\omega(t)$ in $\cal H$. Let us say that the elements $y, y' \in \cal H$ coincide on $\omega$ (and write $y\overset{\omega}=y'$), if there exists $\varepsilon=\varepsilon(\omega,y)>0$ such that $P_{\omega(t)}y=P_{\omega(t)}y'$ for $0\leqslant t< \varepsilon$. Coincidence on $\omega$ obviously is an equivalence relation. The corresponding equivalence class is called the {\it germ} of the element $y$ on the atom $\omega$ and is denoted by $\tilde y(\omega)$. The set of germs ${\cal G}_\omega:=\{\tilde y(\omega)\,|\,\,y \in {\cal H}\}$ forms the {\it stalk} above $\omega$ which obviously has the structure of linear space.

We call the space of ``functions'' $\tilde{\cal H}:=\{\tilde y(\cdot)\,|\,\,y \in {\cal H}\}$ with algebraic operations defined point-wise the {\it model} space and its elements $\tilde y$ {\it models} of elements $y \in \cal H$. Transition to the model is realized by the operator $W: {\cal H} \to\tilde{\cal H}, \,\,W y:=\tilde y(\cdot)$. It is linear and in known applications injective. Non-injectivity of $W$ would mean existence of a non-zero $y\in \cal H$ and a function $\varepsilon=\varepsilon(\omega)$ such that $y \perp \vee_{\omega\in \Omega_{L_0}}\vee_{0 \leqslant t<\varepsilon(\omega)}\omega(t)$, which could be interpreted as absence of {\it completeness} of the system of atoms. In the applications that we know completeness takes place, but whether the same holds in the general case is an open question.

The transition operator $W$ has additional properties, if the space $\tilde{\cal H}$ is equipped with a Hilbert structure. One of the ways to define such a structure is the following. Let ${\rm Ker\,}W=\{0\}$. Take by definition $(\tilde y,\tilde w)_{\tilde {\cal H}}:=(y,w)_{\cal H}$; then $W$ is unitary. If ${\rm Ker\,}W\not=\{0\}$, then by restricting $W$ to ${\cal H}\ominus{\rm Ker\,}W$ we obtain a partial isometry. This trick is used in the model theory (see, e.g., \cite{Shtraus}); it is universal, but not very meaningful. If there was found some canonical Hilbert structure in stalks ${\cal G}_\omega$, then one could hope for realization of $\cal H$ as $\tilde{\cal H}=\oplus\int_{\Omega_{L_0}}{\cal G}_\omega\,d\mu(\omega)$ (with an adequate measure $\mu$) such that $W:{\cal H}\to\tilde{\cal H}$ would be unitary. At present such a structure cannot be seen, but chances appear under additional assumptions about atoms. Studied examples motivate the following heuristic construction.
\smallskip

\noindent{\bf Values on atoms.}
\,\,Recall that in the system (\ref{alpha1})--(\ref{alpha3}) for every control $h \in \cal M$ there exists a corresponding classical solution $u^h \in C^\infty_{\rm loc}\left([0,\infty); {\cal H}\right)$ called smooth wave (see \ref{ssec System alpha}). They constitute reachable sets ${\cal U}_{L_0}^T$ and ${\cal U}_{L_0}$\,\,(see. (\ref{U^t})).

Let the operator $L_0$ be completely non-self-adjoint, so that controllability (\ref{contr}) takes place. Additionally assume that:

\noindent {\bf (A)}\,\,\,there exists a subset $\Omega^e$ of $\Omega_{L_0}$ such that the system of atoms constituting $\Omega^e$ is complete, they all are continuous at zero, and $\omega(0)=\lim \limits_{t \to +0}\omega(t)=\{0\}$;

\noindent {\bf (B)}\,\,\,there exists an element $e \in \cal H$ such that the limits $\lim\limits_{t \to +0}\frac{\|P_{\omega(t)}u\|}{\|P_{\omega(t)} e\|}$ exist and are finite for every $u \in{\cal U}_{L_0},\,\omega \in \Omega^e$. We call such $e$ a {\it gauge} element.

In the general case existence of gauge elements is not proved, however, in examples they can be found, and it is even possible to choose $e \in {\rm Ker\,}L_0^*$. One can call the lineal ${\cal U}_{L_0}$ the {\it smooth structure} determined by $L_0$, owing to the role of the condition {\bf (B)} that we will see below.

Fix $\omega \in \Omega^e$; let $\tilde{\cal U}_{L_0,\omega}:=\{\tilde u(\omega)\,|\,\,u \in {\cal U}_{L_0}\}\subset {\cal G}_\omega$ be the lineal of germs of smooth waves. Define the following sesquilinear form on it
\begin{equation}\label{Hilbert structure}
\langle \tilde u(\omega),\tilde u'(\omega)\rangle\,:=\,\lim\limits_{t \to +0}\frac{(P_{\omega(t)} u, u')}{(P_{\omega(t)} e, e)}\,,\qquad u,u' \in {\cal U}_{L_0},\,\,\, \omega \in \Omega^e\,.
\end{equation}
Consider its (linear) subset $\tilde{\cal U}_{L_0,\omega}^{\,0}:=\{\tilde u(\omega)\,|\,\,u \in {\cal U}_{L_0},\,\,\langle \tilde u(\omega),\tilde u(\omega)\rangle=0\}$ and the factor space ${\cal U}_{L_0,\omega}:=\tilde{\cal U}_{L_0,\omega}/\tilde{\cal U}_{L_0,\omega}^{\,0}$; by $[u](\omega)$ denote the equivalence class of the element $\tilde u(\omega)$. Let us call $[u](\omega)$ the {\it value} of the wave $u\in{\cal U}_{L_0}$ on the atom $\omega$. The form (\ref{Hilbert structure}) induces a natural pre-Hilbert structure on ${\cal U}_{L_0,\omega}$. Taking completion with respect to the corresponding norm, we obtain the Hilbert {\it space of values}. Let us keep the notation ${\cal U}_{L_0,\omega}$ for it. Here we can note that every wave $u\in{\cal U}_{L_0}$ can be represented as $u=u^h(T)$, and so evolution of waves in the system (\ref{alpha1})--(\ref{alpha3}) is reflected in evolution of values $[u^h](\omega, T)$ on atoms $\omega \in \Omega^e$.
\smallskip

\noindent{\bf The wave representation.}
\,\,In addition to {\bf (A)} and {\bf (B)}, let us take another assumption:

\noindent {\bf (C)}\,\,\,There exists a measure $\mu$ on $\Omega_{L_0}$ such that $\mu(\Omega_{L_0}\backslash \Omega^e)=0$ and
\begin{equation}\label{isometry u to u omega}
(u,u')_{\cal H}\,=\,\int_{\Omega_{L_0}}\langle\, [u](\omega), [u'](\omega)\,\rangle\,d\mu(\omega)\,, \qquad u,u' \in{\cal U}_{L_0}\,.
\end{equation}
In the examples we know such measures can be found. It is not known whether the conditions {\bf (A)} and {\bf (B)} guarantee existence of $\mu$ in the general situation.

We call the space ${\cal H}^{\rm w}:=\oplus\int_{\Omega_{L_0}}{\cal U}_{L_0,\omega}\,d\mu(\omega)$ the {\it wave representation} of the original $\cal H$. From the definitions it is clear that the operator $U:{\cal H}\to {\cal H}^{\rm w}$,
$$
{\cal H}\supset{\cal U}_{L_0}\ni u\,\overset{U}\mapsto\,[u](\cdot)\in {\cal H}^{\rm w},
$$
which realizes this representation, is isometric and can be extended to a unitary operator from ${\cal U}_{L_0}$ to the whole $\cal H$. $U$ acts by applying $W$ and factorizing the germs.

The purpose of passing from germs $\tilde u$ to values $[u](\cdot)$ is the following. In all known examples describing elements of $\cal H$ by sections of the bundle $\cup_{\omega \in
\Omega_{L_0}}\{\omega, {\cal G}_\omega\}$ is redundant. Passing to values removes this redundancy, owing to factorization. We show it in the example of the Sturm-Liouville operator.
\medskip

\noindent$\bullet$
\,\,\,Let $L_0$ be the operator (\ref{SL first}). In this case ${\cal H}=L_2({\mathbb R}_+)$.
\smallskip

\noindent{\bf Germs.}\,\,\,Pick an atom $\omega \in\Omega$. Let $y, y' \in {\cal H}$ be two functions. By Theorem \ref{Th 1}, the equality $y\overset{\omega}=y'$ means that $y$ and $y'$ coincide in some neighborhood of the point $x_\omega \in{\bar{\mathbb R}_+}$. Thus the germ $\tilde y(\omega)$ can be canonically identified with the ordinary germ of the function $y(\cdot)$ at the point $x_\omega$, and the model space $\tilde{\cal H}$ with the stalk of square-summable functions above $x_\omega$. This way the stalks ${\cal G}_\omega$ are spaces of infinite dimension.

From the same Theorem \ref{Th 1} it easily follows that the system of atoms composing the wave spectrum $\Omega$ is complete. Thus the operator $W: y \mapsto \tilde y$ is injective. At the same time, owing to ${\rm dim\,}{\cal G}_\omega=\infty$, modeling scalar functions $y$ by the elements of the germ $\tilde y$ is obviously redundant. This motivates passing from germs to values.
\smallskip

\noindent{\bf Values.}
\,\,In our case the condition {\bf (A)} is satisfied.

Let us check the condition {\bf (B)}. The set of smooth waves is ${\cal U}\!=\!C^\infty_{\rm fin}\left({\bar{\mathbb R}_+}\right)$, see (\ref{UT SL}). Recall that the function $\phi$ is the solution of the problem (\ref{phi}). Pick a non-zero element $e \in {\rm Ker\,}L_0^*$. Owing to (\ref{L_0^*}) we have $e=c\phi$ with some constant $c\not=0$. According to the general theory of ordinary differential equations the function $e$ is smooth and can have only simple zeros that can accumulate only to $\infty$. Denote $N^e:=\{x\in {\mathbb R}_+\,|\,\,e(x)=0\}$.

Let $\omega=\{\omega(t)\}_{t\geqslant 0}: \omega(t)=L_2(\{x_\omega\}^t)$ be an atom such that $x_\omega \not\in N^e$. For the smooth waves $u, u' \in \cal U$ we have:
\begin{equation}\label{values SL 1}
\frac{(P_{\omega(t)} u, u')}{(P_{\omega(t)} e,e)}
=\frac{\int_{\{x_\omega\}^t}u(x)\overline{u'(x)}\,dx}
{\int_{\{x_\omega\}^t}e(x)\overline{e(x)}\,dx}
\,\underset{t\to+0}\to\left(\frac{u(x_\omega)}{e(x_\omega)}\right)
\overline{\left(\frac{u'(x_\omega)}{e(x_\omega)}\right)}.
\end{equation}
As we see, the function $e$ is suitable for the role of the gauge element. Taking
$\Omega^e=\Omega_{L_0}\backslash\{\omega\,|\,\,x_\omega \in N^e\}$ we conclude that the condition {\bf (B)} is satisfied.

Owing to (\ref{values SL 1}), we have for the germs $\tilde u(\omega),\tilde u'(\omega)\in\tilde{\cal U}_{\omega} \subset {\cal G}_\omega$
\begin{equation}\label{la ra}
\langle \tilde u(\omega),\tilde u'(\omega)\rangle\,
=\,\left(\frac{u(x_\omega)}{e(x_\omega)}\right)
\overline{\left(\frac{u'(x_\omega)}{e(x_\omega)}\right)}\,.
\end{equation}
Clearly the condition $\tilde u(\omega)\in \tilde{\cal U}_{\omega}^{\,0}$, which by definition means that $\langle\tilde u(\omega),\tilde u(\omega)\rangle=0$, is equivalent to $u(x_\omega)=0$. It also follows from the equality that the correspondence
\begin{equation}\label{Coord values}
{\cal U}_{\omega}
=\tilde{\cal U}_{\omega}/\tilde{\cal U}_{\omega}^{\,0} \ni[u](\omega)\,\,\mapsto\,\, \lim\limits_{t\to+0}\frac{(P_{\omega(t)}u,e)}{(P_{\omega(t)}e,e)}
=\frac{u(x_\omega)}{e(x_\omega)}\in\mathbb C
\end{equation}
is an isometry. Thus for $\omega\in\Omega^e$ the space of values ${\cal U}_{\omega}$ is {\it one-dimensional}. The same correspondence gives the canonical coordinatization of ${\cal U}_\omega$. Other coordinatizations are also possible and have the form $[u](\omega)\mapsto e^{i\theta(\omega)} \frac{u(x_\omega)}{e(x_\omega)}$ with real-valued functions $\theta(\cdot)$.
\smallskip

\noindent{\bf The wave representation.}
\,\,Recall that the measure $\nu$ is defined in Proposition \ref{Prop 3}. For the smooth waves $u,u' \in {\cal U}$ we have:
\begin{align*}
& (u,u')=\int_{{\mathbb R}_+}u(x)\overline{u'(x)}\,dx=\int_{{\mathbb R}_+}
\left(\frac{u(x_\omega)}{e(x_\omega)}\right)
\overline{\left(\frac{u'(x_\omega)}{e(x_\omega)}\right)}\,|e(x)|^2
dx=
\\
&\overset{(\ref{la ra}), (\ref{Coord values})}=\int_{\Omega}\langle
[u](\omega),[u'](\omega)\rangle\,d\mu(\omega)\,,
\end{align*}
where $\mu(\omega):=|e(x_\omega)|^2d\nu(\omega)$. As we see, the condition {\bf (C)} is satisfied and the correspondence
\begin{equation}\label{U map}
{\cal H}=L_2({\bar{\mathbb R}_+})  \ni u\,\,\overset{U}\mapsto\,[u](\cdot) \in L_{2,\,\mu}(\Omega)=:{\cal H}^{\rm w}
\end{equation}
is an isometry. It is defined on smooth waves and can be extended from $\cal U$ to a unitary operator $U: {\cal H}\to {\cal H}^{\rm w}$. The latter gives the wave representation of the elements of the original $\cal H$.
\smallskip

\noindent{\bf The coordinate representation.}
\,\,Coordinatizations of the spectrum (\ref{Coord spectrum}) and of the value spaces (\ref{Coord values}) define an isometry
\begin{equation}\label{V map}
{\cal H}^{\rm w} \ni [u](\cdot)\,\overset{V}\mapsto \,u[\cdot] \in L_{2,\, \rho}({\bar{\mathbb R}_+}[\tau])=:{\cal H}^{\rm c}\,,
\qquad
\end{equation}
where $[u](\tau):=\frac{u(\tau)}{e(\tau)}$ and $d\rho:=|e(\tau)|^2{d\tau}$. It gives the wave representation of elements of space ${\cal H}^{\rm w}$.

Composition $Y:=VU$, $Y:{\cal H}\to{\cal H}^{\rm c}$ can be extended from smooth waves to a unitary operator which maps functions from the original $L_2\left({\bar{\mathbb R}_+}[x]\right)$ to functions from $L_{2,\, \rho}({\bar{\mathbb R}_+}[\tau])$ by the rule
\begin{equation}\label{Y}
\left(Yy\right)(\tau)\,=\,y[\tau]:=\frac{y(\tau)}{e(\tau)}\,,\qquad \tau \geqslant 0\,.
\end{equation}

Obvious similarity of the original space to the space of coordinate repre\-sen\-tation and simplicity of the correspondence $y\mapsto y[\cdot]$ are important facts used for solving inverse problems.

\subsection{The operator $\tilde L_0^*$}\label{Operator tilde L0*}

Let us return to the general case that was considered in the beginning of Section \ref{ssec Space tilde cal H}. Assume that controllability $\overline{\cal U}_{L_0}=\cal H$ takes place and the operator of transition to the model $W$ is injective. In this situation the set of pairs $\left\{\{Wu,WL_0^*u\}\,|\,\,u \in {\cal U}_{L_0}\right\}$ is a graph of an operator which acts in the model space $\tilde {\cal H}$. We call it the {\it wave model} of the operator $L_0^*$ and denote by $\tilde L_0^+$. Note that it would be more consistent to talk about the model not of $L_0^*$ itself, but of its wave part $L_0^*|_{{\cal U}_{L_0}}$ (see the remark in Section \ref{ssec System alpha}). We ignore this inaccuracy in order to not overload terminology.

Since every smooth wave is $u=u^h(T)$ and
$$
L_0^*u^h(T)\overset{(\ref{alpha1})}=-u^h_{tt}(T)\overset{(\ref{weak solution u^f})}=-u^{h_{tt}}(T)\,,
$$
the wave model can be defined as operator with the graph
\begin{equation*}
{\rm graph\,}\tilde
L_0^+\,=\,\left\{\{Wu^h(T),-Wu^{h_{tt}}(T)\}\,|\,\,h \in {\cal
M},\,\,T\geqslant 0\right\}\,.
\end{equation*}

It is hard to expect for rich properties of the model in such generality. The conjecture is its {\it locality}: if $\tilde y \in {\rm Dom\,}\tilde L_0^+$ and $\tilde y|_{A}=0$ on an open set ${A}\subset \Omega_{L_0}$, then $\tilde L_0^+\tilde y|_{A}=0$. In known examples locality takes place.

If the wave representation (\ref{U map}) is defined, then the corresponding version of the wave model appears, $L_{0\, \rm w}^+:= U L_0^*|_{{\cal U}_{L_0}} U^*$, which acts in the space ${\cal H}^{\rm w}$. It is defined by its graph
\begin{equation*}
{\rm graph\,} L_{0\, \rm w}^+\,=\,\left\{\{Uu^h(T),-Uu^{h_{tt}}(T)\}\,|\,\,h \in {\cal M},\,\,T\geqslant 0\right\}\,.
\end{equation*}

\noindent$\bullet$\,\,In the case of the Sturm-Liouville operator the coordinate realization of the wave model $L_{0\,\rm c}^+:= Y L_0^*|_{{\cal U}_{L_0}} Y^*$ is defined correctly. It acts in the space ${\cal H}^{\rm c}=L_{2, \,\rho}({\bar{\mathbb R}_+}[\tau])$, is defined by its graph
\begin{equation*}
{\rm graph\,} L_{0\, \rm c}^+\,=\,\left\{\{Yu^h(T),-Yu^{h_{tt}}(T)\}\,|\,\,h \in {\cal M},\,\,T\geqslant 0\right\}
\end{equation*}
and, owing to (\ref{L_0^*}) and (\ref{Y}), is the differential operator
\begin{align}
\notag & \left(L_{0\, \rm c}^+y[\cdot]\right)[\tau]\,=\,\left\{\frac{1}{e(\tau)}\,\left(-\frac{d^2}{d\tau^2}
+q(\tau)\right)e(\tau)\right\}y[\tau]\,=
\\
\label{L0*coord} & =
-\,y''[\tau]+p(\tau)\,y'[\tau]+Q(\tau)y[\tau]\,,\qquad \tau\geqslant 0,\,\,\tau\in\Omega^e\,.
\end{align}
with the coefficients
\begin{align}\label{p and Q}
p(\tau):=-2\,\frac{e'(\tau)}{e(\tau)}\,, \quad Q(\tau):=q(\tau)-\frac{e''(\tau)}{e(\tau)}\,.
\end{align}
This operator is not closed, but its closure $L_{0\,\rm c}^*=\overline{L_{0\, \rm c}^+}$ is unitarily equivalent to the operator $L_0^*$ by the remark at the end of Section \ref{ssec Controllability}. The construction of the canonical Green's system from Section \ref{ssec System Gr_L_0} results in the following (recall that $(Yy)(\tau)=\frac{y(\tau)}{e(\tau)},\ (Y^*y)(\tau)=y(\tau)e(\tau)$):
\begin{equation}
\begin{array}{l}
K_c=YK=\{\text{const}\},
\\
\Gamma_{1\,\rm c}=Y\Gamma_1Y^*: y\mapsto-y(0),
\\
\Gamma_{2\,\rm c}=Y\Gamma_2Y^*: y\mapsto\frac{y'(0)}{\eta'(0)}.
\end{array}
\end{equation}

\subsection{The inverse problem}\label{Inverse Problem}
The functional model $\tilde L_0^+$ gives a unified approach to a rather wide class of boundary inverse problems. Putting off generalizations, we demonstrate the idea of the approach with our example.
\smallskip

\noindent$\bullet$\,\,{\bf Auxiliary model.}\,\,Consider the boundary problem
\begin{align}
\label{Psi 1} & - \psi'' + q\psi\,=\,\lambda \psi, && x>0
\\
\label{Psi 2} & \psi(0)=0, \quad \psi'(0)=1\,.
\end{align}
Its solution $\psi=\psi(x,\lambda)$ is a function that is smooth in $x$ and entire in $\lambda \in \mathbb C$. In particular, for $q=0$ we have $\psi=\frac{\sin\sqrt{\lambda}\,x}{\sqrt{\lambda}}$.

Our operator $L_0$ has defect indices $(1,1)$. Therefore there exists a unique {\it spectral function} $\sigma$ such that the formulas
$$
\check y(\lambda)=\int_0^\infty y(x)\,\psi(x,\lambda)\,dx, \qquad y(x)=\int_{-\infty}^\infty \check y(\lambda)\,\psi(x,\lambda)\,d\sigma(\lambda)
$$
establish an isometry of the spaces $L_2({\mathbb R}_+[x])$ and $L_{2,\,\sigma}({\mathbb R}[\lambda])$ (cf. \cite{Naimark}, Chapter VIII). In other words, the operator $\Phi: y \mapsto \check y$ is unitary.

Let us find $\Phi$-representation of waves from the problem (\ref{alpha1*})--(\ref{alpha3*}) with a smooth control $f \in\dot{\cal M}$. For $\check u^f(\cdot, t)=\Phi u^f(\cdot,t)$, using finiteness of the support of $ u^f(\cdot, t)$, we have:
\begin{align*}
& \check u^f_{tt}(\lambda,t)=\int_0^\infty
u^f_{tt}(x,t)\psi(x,\lambda)\,dx\overset{(\ref{alpha1*})}=\int_0^\infty
\left[u^f_{xx}(x,t)-q(x)u^f(x,t)\right]\psi(x,\lambda)\,dx\\
& =  -u^f_x(0,t)\psi(0,\lambda)+u^f(0,t)\psi'(0,\lambda)+
\int_0^\infty u^f(x,t)\left[\psi''(x,\lambda)-q(x)\psi(x,\lambda)\right]\,dx
\\
& \overset{(\ref{alpha1*}), (\ref{alpha3*}), (\ref{Psi 1}),
(\ref{Psi 2})}=f(t)-\lambda \int_0^\infty
u^f(x,t)\psi(x,\lambda)\,dx = f(t)-\lambda\check u^f(\lambda,t)\,.
\end{align*}
Integrating and using (\ref{alpha2*}) we find:
\begin{equation}\label{u check 1}
\check
u^f(\lambda,t)\,=\,\int_0^t\lambda^{-\frac{1}{2}}\,
\sin[\lambda^{\frac{1}{2}}\,(t-s)]f(s)\,ds\,,\qquad t \geqslant 0\,.
\end{equation}

For the inverse problem the operator $\check L_0^+:=\Phi\, L_{0}^*|_{\cal U}\,\Phi^*$  plays the role of an {\it auxiliary model} of the original $L_0^*$. It acts in $L_{2,\,\sigma}({\mathbb R}[\lambda])$ by the rule
\begin{align}
\notag & \check L_0^+\check
u^f(\lambda,t)\,=\left[L_0^*u^f(\cdot,t)\right]^{\check{}}(\lambda)=-[u^f_{tt}(\cdot,t)]^{\check{}}(\lambda)=
-[u^{f_{tt}}(\cdot,t)]^{\check{}}(\lambda)\,=
\\
\label{u check 2} &=-\check u^{f_{tt}}(\lambda,t)\overset{(\ref{u check 1})}
=-\,\int_0^t\lambda^{-\frac{1}{2}}\,\sin[\lambda^{\frac{1}{2}}\,(t-s)]f''(s)\,ds
\end{align}
and is defined by its graph
\begin{equation}\label{graph check L0*}
{\rm graph\,} \check L_0^+\,=\,\left\{\{\check u^f(\lambda,T),-\check u^{f_{tt}}(\lambda,T)\}\,|\,\,f \in \dot{\cal M},\,\,T\geqslant 0\right\}\,.
\end{equation}

\noindent{\bf Recovering the potential.}
\,\,The classical inverse spectral problem for the Sturm-Liouville operator on the half-line is to determine the potential $q|_{x \in {\bar{\mathbb R}_+}}$ from the given spectral function $\sigma|_{\lambda \in\mathbb R}$ (cf. \cite{Naimark}, Chapter VIII). We can solve this problem by the following scheme.

\noindent{\bf Step 1.}\,\,Using (\ref{u check 1}) and (\ref{u check 2}), find the operator $\check L_0^+$ in the space $L_{2,\,\sigma}({\mathbb R}[\lambda])$ from its graph
(\ref{graph check L0*}).

\noindent{\bf Step 2.}\,\,Construct the wave model of the operator $\check L_0^+$ and consider its coordinate realization. Owing to invariance of the construction of the wave model this leads to the operator $L_{0\,\rm c}^*$ acting in $L_{2,\,\rho}({\mathbb R}_+[\tau])$.

\noindent{\bf Step 3.}\,\,From the representation (\ref{L0*coord}) find the coefficients $p$ and $Q$. From (\ref{p and Q}) find $e(\tau)=C{\rm exp\,}\{-\int_0^\tau \frac{p(s)}{2}\,ds\}$.
Finally, recover $q(\tau)=Q(\tau)+\frac{e''(\tau)}{e(\tau)}$.
\smallskip

Owing to complexity of the construction of the wave model, this scheme is, of course, too involved compared with the classical procedure that uses the Gelfand-Levitan equation \cite{Naimark}.
In the present paper we only want to demonstrate, on a relatively simple example of the Sturm-Liouville operator, the construction in all the details and to show how it solves inverse problems.

At the same time, the wave model has some advantages. It can be used for solving problems with {\it any} data, if they only determine the operator $L_0$ (or, equivalently, $L_0^*$) up to unitary equivalence. Spectral data, scattering data, Weyl function, characteristic function (cf. \cite{Shtraus}, \cite{MMM}, \cite{Ryzh}) can be considered as such. Universality of the model makes it unnecessary to convert the data of one type into another. Besides that, the wave model is efficient for recovering objects of greater complexity, Riemmanian manifolds \cite{JOT}.

In the future we plan to study the construction of the wave model itself, as well as its possible applications. It would be interesting to construct, on the basis of the wave model, a functional model of an abstract Green's system ${\mathfrak G}_{L_0}$, see Section \ref{ssec System Gr_L_0}. This interest is motivated by existence of the {\it boundary} of the wave spectrum (\ref{boundary}).

There are relations, which are worth being studied, between the wave model and operator $C^*$-algebras \cite{JGPh}. The source of these relations is the correspondence $\omega \leftrightarrow \tau_\omega$ between atoms and eikonals, see (\ref{eikonal}).

\end{document}